\pdfoutput=1
\documentclass[10pt]{article}
\usepackage{graphicx}
\usepackage{amssymb}
\usepackage{epstopdf}
\usepackage{amsmath}
\usepackage{amsfonts}

\def\phi{{\varphi}}

\DeclareSymbolFont{AMSb}{U}{msb}{m}{n}
\DeclareMathSymbol{\N}{\mathbin}{AMSb}{"4E}
\DeclareMathSymbol{\Z}{\mathbin}{AMSb}{"5A}
\DeclareMathSymbol{\R}{\mathbin}{AMSb}{"52}
\DeclareMathSymbol{\Q}{\mathbin}{AMSb}{"51}
\DeclareMathSymbol{\I}{\mathbin}{AMSb}{"49}
\DeclareMathSymbol{\C}{\mathbin}{AMSb}{"43}

\def\be{\begin{equation}}
\def\ee{\end{equation}}
\def\ber{\begin{eqnarray}}
\def\eer{\end{eqnarray}}

\def\beq{\begin{equation}}
\def\eeq{\end{equation}}

\begin{document}

\addtolength{\textheight}{0 cm} \addtolength{\hoffset}{0 cm}
\addtolength{\textwidth}{0 cm} \addtolength{\voffset}{0 cm}

\newenvironment{acknowledgement}{\noindent\textbf{Acknowledgement.}\em}{}

\setcounter{secnumdepth}{5}
 \newtheorem{proposition}{Proposition}[section]
\newtheorem{theorem}{Theorem}[section]
\newtheorem{lemma}[theorem]{Lemma}
\newtheorem{coro}[theorem]{Corollary}
\newtheorem{remark}[theorem]{Remark}
\newtheorem{claim}[theorem]{Claim}
\newtheorem{conj}[theorem]{Conjecture}
\newtheorem{definition}[theorem]{Definition}
\newtheorem{application}{Application}

\newtheorem{corollary}[theorem]{Corollary}

\title{Metric Selfduality and Monotone Vector Fields on Manifolds}
\author{Nassif  Ghoussoub\\
{\it\small Department of Mathematics}\\
{\it\small  University of British Columbia}\\
{\it\small Vancouver BC Canada V6T 1Z2}\\
{\it\small nassif@math.ubc.ca}\vspace{1mm}\and
Abbas Moameni\footnote{Both authors were partially supported by grants
from the Natural Sciences and Engineering Research Council of Canada.}
%\thanks{Partially supported by a grant from the Natural Sciences and Engineering Research Council of Canada.}
\hspace{2mm}\\
{\it\small School of Mathematics and Statistics}\\
{\it\small Carleton University}\\
{\it\small Ottawa, ON, Canada K1S 5B6 }\\
{\it\small momeni@math.carleton.ca}\\\\
}
\maketitle

%
%\date{Revised November 30, 2011}
%\today

\begin{abstract} We develop a ``metrically selfdual" variational calculus for  $c$-monotone vector fields between general manifolds $X$ and $Y$, where $c$ is a coupling on $X\times Y$.  
Remarkably, many of the key properties of classical monotone operators known to hold in a linear context,  extend to this non-linear setting. This includes an integral representation of $c$-monotone  vector fields in terms of $c$-convex  selfdual Lagrangians, their characterization as a partial $c$-gradients of antisymmetric Hamiltonians, as well as the property that these vector fields are generically single-valued. 
We also use a symmetric Monge-Kantorovich transport to associate to any measurable map  its closest possible $c$-monotone ``rearrangement". 
We also explore how this metrically selfdual representation can lead to a global variational approach to the problem of inverting $c$-monotone maps, an approach that has proved efficient for resolving non-linear equations and evolutions driven by monotone vector fields in a Hilbertian setting. 
 
\end{abstract}
\tableofcontents
\section{Introduction and main results} 
 
Many aspects of convexity theory such as Fenchel-Legendre duality, subdifferentiability, and cyclic monotonicity have been extended to settings where the usual linear duality $\langle x, x^*\rangle$ between a Banach space $X$ and its dual $X^*$ is replaced by a general coupling $c(x, y)$ of two arbitrary sets $X$ and $Y$. These nonlinear ``metric" generalizations of convexity and cyclic monotonicity were mostly motivated by problems in Riemannian geometry \cite{Mn},  mathematical economics \cite{Ca}, \cite{ML}, and by the Monge-Kantorovich theory of mass transport corresponding to general cost functions \cite{V1}. For example, McCann's extension of Brenier's theorem \cite{Br} to manifolds required that the scalar product in the linear theory be replaced by $c(x,y)=-d^2(x,y)$, where $d$ is the Riemannian metric and where convexity is replaced by the concept of {\it ``$c$-convexity"}  described below. What is remarkable is that many of the key structural results known to hold under assumptions of classical convexity and cyclic monotonicity on Euclidean space extend to this metric setting. This has had major impact on non-linear analysis, differential geometry and their applications. See for example the books of Villani \cite{V1, V2}). In this paper, we show that similar metric extensions hold for the notions of selfduality and for vector fields that are merely 2-monotone. 

On the other hand, the most natural extensions of gradient flows of convex energies on a Hilbert space to curved manifolds seem to be those corresponding to  {\it ``pathwise convex"} and not necessarily {\it ``metrically convex"} energies. Indeed, Otto's calculus \cite{V1} allows the rewrite of several non-linear evolution equations as gradient flows of geodesically convex free energy functionals on the Wasserstein manifold. Unfortunately, the global variational methods that characterize the success of convex analysis on linear spaces do not readily extend to non-linear settings. Instead, time discretization methods had to be used to circumvent the lack of global variational principles. See for example  the penetrating study of Ambrosio-Gigli-Savar\'e \cite{AGS}, or the more recent comprehensive notes of Ambrosio-Gigli \cite{AG}. Our quest for global variational methods in nonlinear settings eventually led us to interesting connections between these two notions of convexity, 
including a new criterion to deduce the metric convexity of pathwise convex functionals, which  normally is not easy to verify.
 
We begin by stating the main properties of (sub)-gradients of convex functions and monotone vector fields that we will be extending to a non-linear setting. First, we recall two fundamental notions of monotonicity.\\
 A --possibly set valued-- map  $T: Dom (T)\subset X \to X^*$  
 is said to be 
{\it cyclically monotone} (or {\it $n$-cyclic monotononicity for every $n$}) if its graph $G(T)=\{ (x,p)\in X\times X^*; p\in Tx\}$ satisfies the property that 
 for any finite number of points $(x_i, p_i)_{i=0}^n$  on $G(T)$ with $x_0=x_n$, we have
\begin{equation}
\hbox{$\sum\limits_{i=1}^n\langle p_i, x_i-x_{i-1} \rangle \geq 0$.}
\end{equation}
On the other hand, $T$ is said to be {\it monotone} if it is only {\it 2-cyclically monotone}, meaning that its graph satisfies \begin{equation}
\hbox{$\langle x_1-x_2, p_1 - p_2 \rangle \geq 0$ for every $(x_1,p_1)$ and $(x_2, p_2)$ in $G(T)$.}
\end{equation}
$T$ (or its graph $G(T)$) is said to be {\it maximal cyclically monotone} (resp., {\it maximal monotone}) in $X\times X^*$, if it has no proper cyclically monotone  (resp., monotone) extension in  $X \times X^*.$ 

Here are the structural results enjoyed by such vector fields that we plan to extend to a non-linear setting.\\

\noindent {\bf A. Main properties of cyclically monotone maps:}  
\begin{enumerate}
\item {\it Integral representation} (Rockafellar \cite{Rock}): $T$ is a (maximal) cyclically monotone vector field if and only if there exists a convex lower semi-continuous function $\phi$ such that $Tx=\partial \phi (x)$ for all $x$ in its domain, 
where $\partial \phi$ is the subdifferential of $\phi$. 
\item  {\it Range of cyclically monotone maps} (Asplund \cite{Ph}): Proper lower semi-continuous convex functions on a reflexive space are Fr\'echet-differentiable on a dense $G_\delta$ subset of their domains. 

\item {\it Optimal mass transport} (Brenier \cite{Br}, Gangbo \cite{Gg}):  For any continuous measure $\mu$ and any non $\mu$-degenerate measurable map  $T: \Omega \to \R^d $, there is a cyclically monotone vector field $T_{\infty}: \Omega \to \R^d$ of the form $T_{\infty}=\partial \psi$, where $\psi$ minimizes  the functional $I(\phi) =\int_\Omega (\phi (x) +\phi^*(Tx))\, d\mu$ over all convex functions on $\R^d$. Here $\phi^*$ is the Fenchel-Legendre dual of $\phi$.

\item {\it Variational inversion of a gradient field}: If $\phi$ is a convex function associated to a cyclically monotone map $T$, then under mild coercivity conditions, the infimum of $I_p(u)=\phi (u)-\langle u,p\rangle$ is attained at a point $\bar u$ that solves the equation $p\in T{\bar u}$.

\end{enumerate}

\noindent The theory of monotone operators dates back to the sixties \cite{Brezis}, but the existence of a ``convex integral representation" for a monotone map $T$ is more recent and was started by Fitzpatrick \cite{F}. The related notion of a {\it selfdual Lagrangian} was introduced independently in \cite{GT1} in order to address a conjecture by Brezis-Ekeland regarding a variational characterization of gradient flows \cite{BE}. That one can refine the Fitzpatrick representation to obtain a representation by a selfdual Lagrangian was discovered by several authors (see for example the book \cite{G10}).  The theory was developed further in \cite{G1, GM2, GM3, GM4} in order to deal with other PDEs and evolution equations. \\
We recall that a {\it selfdual Lagrangian} is a 
lower semi-continuous convex  function $L:X\times X^*\to \R\cup \{+\infty\}$  which satisfy the following {\it selfduality} conditions:
\begin{equation}
\hbox{$L^*( p, u) = L(u, p )$   for all $(p,u)\in X^{*}\times X$.}
\end{equation}
Here $X$ is a reflexive Banach space and $L^*$ is the Legendre transform of $L$ in both variables, that is
\[
L^*(p,u)=\sup\left\{ \langle p,y\rangle +\langle u,q\rangle -L(y, q);\, (y, q)\in X\times X^*\right\}.             
\]
It is easy to see that such Lagrangians satisfy 
$
L(u, p) -\langle u, p\rangle \geq 0$ for all $(u,p)\in X\times X^*$, 
and 
that $L(u,p)-\langle u,p\rangle=0$ if and only if $(p,u)\in \partial L(u,p),$ where $\partial L$ is the subdifferential of $L$ in both variables.
 We can therefore associate  to $L$ the  following --possibly set valued-- vector field,
\begin{equation} \label{selfdual.vector}
\hbox{$u\to {\bar\partial} L (u):=\{p\in X^*; \, L(u,p)-\langle u,p\rangle =0\}=\{p\in X^*; \, (p,u)\in \partial L(u,p)\}$.}
\end{equation}
 We can now state the analogous results for $2$-monotone maps. For details, we refer to the book \cite{G10} and the references therein. \\

\noindent {\bf B. Main properties of monotone maps:}  \begin{enumerate}
\item {\it Convex representation of monotone fields} (Fitzpatrick \cite{F}, Ghoussoub \cite{GBrowder}):  $T: Dom (T)\subset X\to 2^{X^*}\setminus \{\emptyset\}$ is a maximal monotone map if and only if there exists a selfdual Lagrangian $L: X\times X^* \to \R$ such that $T={\bar\partial} L.$

\item  {\it Range of monotone maps} (Kenderov \cite{Ph}): Maximal monotone maps on a reflexive Banach space are single-valued on a dense $G_\delta$ subset of their domains. 

\item {\it Symmetric optimal mass transport} (Ghoussoub-Moameni \cite{Gh-Mo1}): For any continuous measure $\mu$ and any non $\mu$-degenerate map $T: \Omega \to \R^d $, there is a monotone vector field $T_2: \Omega \to \R^d$ of the form ${T_2}=\bar \partial M$, where $M$ minimizes  the functional $I(L) =\int_\Omega L (x, Tx) \, d\mu$ over all selfdual Lagrangians on $\R^d \times \R^d$.

\item {\it Variational inversion of a monotone field} (Ghoussoub-Tzou \cite{GT1, G1}): If $L$ is a selfdual Lagrangian associated to a monotone map $T$ on a reflexive Banach space $E$, then under mild coercivity conditions, the infimum of $I_p(u)=L(u,p)-\langle u,p\rangle$ over $E$ is zero and is attained at a point $\bar u$ that solves the equation $p\in T{\bar u}$. 

\end{enumerate}

Our goal here is to check to what extent the above results extend to a nonlinear setting, that is when the scalar product $\langle x, y\rangle$ on phase space $X\times X^*$ is replaced by a general ``coupling" or ``cost function" $c(x,y)$ between more general manifold products $X\times Y$. We shall see that ``metric convexity"  and ``metric monotonicity"  seem to be the right context for the extension of B-1, B-2, B-3. However, the non-linear extension of the variational principle B-4 requires an assumption of "arc-wise convexity." See Section 6. 
 
Here are some of the notions under study. Various aspects of $c$-monotonicity with respect to a coupling $c:X\times Y \to \R$ between  two arbitrary sets $X$ and $Y$ have been already introduced in several contexts. We refer to \cite{Ca, ML, RR} for details.
We recall that a subset $M$ of $X \times Y$ is said
to be 
\begin{itemize}
\item {\it $c$-cyclically monotone of order $n$}, if for any set
of pairs $\{(u_i, v_i)\}_{i=1}^n \subseteq M$ with 
$ u_{n+1} = u_1$, we have 
\begin{equation}\sum_{i=1}^n
[c(u_i, v_i) - c(u_{i+1}, v_i)] \geq 0.
\end{equation}
 
\item  {\it $c$-monotone} if it is $c-$cyclically monotone of order $2$, i.e., if  
\begin{equation}
c\big (u_1,v_2\big)+c\big ( u_2,v_1\big)\leq c\big (u_1,v_1\big)+c\big ( u_2,v_2\big), \quad \forall (u_1, v_1), (u_2, v_2) \in M.
\end{equation}

\item $M$ is said to be {\it maximal $c-$cyclically monotone} (resp., {\it maximal $c-$monotone}) in $X\times Y$, if it has no proper $c-$cyclically
monotone  (resp., $c-$monotone) extension in $X \times Y.$ 

\item A set-valued map $T: {\rm Dom}(T)\subset X \to 2^Y\setminus \{\emptyset\}$ is said to be {\it maximal $c-$monotone} (resp., {\it maximal $c-$cyclically monotone}) if its graph $M=G(T)$ is a maximal $c-$monotone (resp., {\it a maximal $c-$cyclically monotone}) subset of $X \times Y.$
\end{itemize}
Note that if  $\langle\, .\, ,.\rangle$ is an inner product on $X \times X$, where $X$ is a Hilbert space, and if one consider the coupling $c(x,y)=-d^2(x,y)/2$, where $d$ is the metric on $X \times X$ induced by this inner product, i.e., $d(x,y)^2=\langle x-y, x-y\rangle$, then it is easy to see that a $c$-monotone map $u$ is necessarily monotone in the classical sense. 

As to the ``metric extension" of the notion of convexity, it goes as follows: Let $C: U\times V \to \R$ be a coupling between two arbitrary spaces $U$ and $V$. For $f: U \to \R \cup \{+ \infty\}$, one can define the $C$-conjugate of $f$ by
\[f^C: V \to \R \cup \{+ \infty\}, \qquad f^C(v)=\sup_{u \in U} \{C(u,v)-f(u)\},\]
and its doubly $c$-conjugate as, 
\[f^{CC}: U \to \R \cup \{+ \infty\}, \qquad f^{CC}(u)=\sup_{v \in V} \{C(u,v)-f^c(v)\}.
\]
A function $f$ is then said to be {\it $C$-convex} if it is equal to its double $C$-conjugate. The $C$-subdifferential of $f$ is the set-valued  map  $\partial_C f: U \to 2^V,$  defined for any $u_0\in U$, by 
\[
\hbox{$\partial_Cf(u_0)=\{v_0 \in V; f(u_0)-f(u) \leq C(u_0,v_0)-C(u,v_0)$ for all $u \in U\}.$
}\]
Just like in the case of a linear coupling, it is easy to see that the  $C$-subdifferentials of $C-$convex functions are maximal $C$-cyclically monotone. Conversely, the same proof as Rockafellar's in the linear setting \cite{Rock} gives that if $M$ is a maximal $C-$cyclically monotone subset of $U \times V$,  then there exists a $C-$convex function $\phi: U \to \R \cup \{+ \infty\}$ such that $M= {\rm Graph} \big (\partial_C\phi\big)$. Properties A-2 and A-3 can also be extended to $c$-cyclically monotone maps (see \cite{V1}), while the analogue of  A-4 will be discussed in Section 6.

One of the objectives of this paper is to define a non-linear version of selfduality that could still provide integral representations for $c$-monotone sets.  The ultimate goal is to couple this representation with a variational principle analogous to B-4, that will allow for resolving equations driven by $c$-monotone vector fields. Here is the non-linear version of selfduality that we propose. 

Let $X$ and $Y$ be two arbitrary sets and  $c: X \times Y \to \R$ be  a coupling.  We consider a new coupling $C$ on the symmetrized space $U\times V$, where $U:=X \times Y$ and $V:=Y \times X$, via the formula
\begin{equation}
C\big((x_1,y_1), (y_2,x_2)\big)=c(x_1,y_2)+c(x_2,y_1). 
\end{equation} 
Say that a function $L: X\times Y \to \R \cup \{+ \infty\}$ is {\it $C$-selfdual} if
\begin{equation}
\hbox{$L^C(y,x)=L(x,y)$  for all $(x,y)\in X \times Y$,}
\end{equation}
 where $L^C: V= Y \times X \to \R \cup \{+ \infty\}$ is the $C$-conjugate of $L$ defined for $v=(y,x) \in V=Y \times X$ by
\begin{eqnarray} \label{kdp}
L^C(v)&=&\sup \{ C(u,v) -L(u);  u \in U=X \times Y\}\\
&=&\sup \{ C((x_1,y_1), (y,x))-L(x_1,y_1); (x_1,y_1) \in X \times Y\}. 
 \end{eqnarray}
 It is easy to see that if $L$ is a  {\it $C$-selfdual} Lagrangian, then we have 
 \begin{equation}
 \hbox{$L(x,y)\geq c(x,y)$ for all $(x,y)\in X\times Y$,}
 \end{equation}
 and that 
 \begin{equation}
 \hbox{$L(x,y) = c(x,y)$  if and only if $(y,x)\in \partial_CL(x,y)$.}
 \end{equation}
We are interested in the possibly set-valued map $\bar \partial_c L: X \to 2^Y\setminus \{\emptyset\}$ defined by
\begin{equation}
\bar \partial_c L(x)=\{y \in Y;\, \, L(x,y)=c(x,y)\}
\end{equation}
and its domain $D_{c,L}$ consisting of all $x \in X$ such that $\bar \partial_c L(x) \not=\emptyset.$  In other words,
\[D_{c,L}=\{x \in X; \, \exists  \, y\in Y \text{ such that } L(x,y)=c(x,y) \}.
\]
We shall say that a function $H$ on $X\times X$ is a {\it sub-antisymmetric Hamiltonian} $H$ if $H(x,x)=0$ and $H(x,y)+H(y,x)\leq 0$ for all $(x,y)$ in the domain of $H.$

 Remarkably, the first three results reminiscent of the linear theory, i.e., B-1, B-2, B-3 extend to this setting, starting with the following selfdual representation that will be established in Section 3.    

  \begin{theorem}\label{three} Let $X$ and $Y$ be two sets and $c: X \times Y \to \R$ be a coupling. The following assertions are then equivalent:
  \begin{enumerate}
  \item $T$ is a maximal $c-$monotone map from  ${\rm Dom}(T)\subset X$ to $2^Y\setminus \{\emptyset\}$.
  
  \item There exists a $C-$selfdual function $L: X \times Y \to \R \cup \{+ \infty\}$ such that $T=\bar \partial_c L$ on ${\rm Dom}(T)$.
  
  \item  There exists an sub-antisymmetric Hamiltonian $H$ on $X\times X$ that is $c$-convex in the second variable such that $Tx = \partial^c_2 H(x,x)$ for all $ x \in {\rm Dom}(T)$. 
\end{enumerate}
\end{theorem}
In Section 4, we study cases where the range of a $c$-monotone map is single-valued. 
In view of the above representation, 
this property is directly linked to the differentiability of antisymmetric functions on  smooth  manifolds. Here we need to assume that the coupling $c$ satisfies the following properties: 
\begin{itemize}
\item {\em The twist condition}, i.e.,  
\begin{equation}\label{cost1}
\hbox{ $D_1 c(x,y_1)=D_1 c(x,y_2)$ implies that $y_1=y_2.$}
\end {equation}
\item For each measurable map $f: X \to Y,$ there exists $p>1$ and a  function  $\eta \in L_{loc}^p(X)$ such that 
\begin{equation} \label{cost2}
|D_1c\big (x,f(x)\big )| \leq |\eta(x)|, \qquad \forall x\in X.
\end{equation}
\end{itemize}
 \begin{theorem}\label{sing.0}
 Let $X$ be  a second countable $C^1$ manifold of dimension $d$ equipped with its volume measure $\mu$, and let $Y$ be a polish space. Assume that $c: X \times Y \to \R$ is a  measurable coupling that is differentiable with respect to the first variable and satisfying (\ref{cost1}) and (\ref{cost2}). 
 Then, any graph measurable  $c-$monotone map $T: Dom(T) \subset  X \to 2^Y\setminus \{\emptyset\}$ is single-valued on its domain  up to a $\mu$-null set. 
  \end{theorem}
An immediate corollary is the following extension of a result by Champion-DePascale \cite{de-ch}, who showed that under rather  similar conditions on the cost function $c$,  it suffices that a transport plan be only $c$-monotone (and not necessarily $c$-cyclically monotone) to insure that it is supported on a graph of a Borel map. The following corollary strengthens that result by showing that the graph is actually a partial gradient of an anti-symmetric Hamiltonian. 
 
\begin{corollary}\label{champion} Let $X$ be a second countable $C^1$ manifold of dimension $d$ and let $Y$ be a
Polish space. Suppose $c$ is a cost function on $X\times Y$ that is differentiable with respect to the first  variable and satisfying (\ref{cost1}) and (\ref{cost2}). Let $\mu$ (resp., $\nu$) be Borel probabilities on $X$ (resp., $Y$) such that $\mu$ is absolutely continuous with respect to the volume measure on $X$. Let $\gamma \in \Gamma(\mu, \nu)$ be a transport plan (i.e., a probability measure on $X\times Y$ with marginals $\mu$ and $\nu$) that is concentrated on a Borel measurable  $c-$monotone subset of $X\times Y$, then
\begin{enumerate}
\item $\gamma$ is necessarily concentrated on the graph of a Borel function $T$ from $X$ to $Y$. 
\item There exists an anti-symmetric Hamiltonian $H$ on $X\times X$, that is $c$-convex in the second variable such that $Tx=\partial^c_2 H(x, x)$ for $\mu$ almost all $x\in X$.
\item If $H$ is locally Lipschitz, and its $c$-conjugate with respect to the second variable is continuous, then for  $\mu$ almost all $x\in X$, $H$ is differentiable with respect to the second variable and 
\begin{equation}
\nabla_2H(x, x)=\nabla_1c(x, Tx).
\end{equation}
\end{enumerate}
\end{corollary} 
In the case when $X$ is a Riemannian manifold equipped with the coupling  $c(x,y)=-d^2(x,y)/2$ induced by its metric $d$, McCann \cite{Mn} had shown that under suitable conditions on the manifold, a continuous map $T: X \to X$ is $c$-cyclically monotone if and only if it can be written as $Tx=\exp_x[\nabla \phi(x)]$, where $\phi: X \to \R$ is a differentiable $c$-convex function.
 One of the applications of our results is the following characterization of $c$-monotone maps on manifolds.

\begin{corollary} Let $(M,g)$ be a connected compact $C^3$-smooth Riemannian manifold without boundary, equipped with a Riemannian distance $d(x,y)$  and its volume measure $\mu$.  Set $c(x,y)=-d^2(x,y)/2$ and consider $T: X \to X$ to be a continuous  map. Then, $T$ is $c$-monotone if and only if there exists an sub-antisymmetric Hamiltonian $H$ on $X\times X$, that is $c$-convex in the second variable such that  $Tx=\exp_x[\nabla_2 H(x,x)]$ for  $\mu$-almost 
every   $x \in X.$
\end{corollary}

In Section 5, we use a symmetric  version of Monge-Kantorovich theory to associate to any vector field $u$ a $c$-monotone map. We do that by considering the class ${\mathcal L}$  of $C$-selfdual functions on $X \times Y,$ i.e.
\[{\mathcal L}=\big \{L: X \times Y \to \R \cup \{+ \infty\}; \, L^C(y,x)=L(x,y), \quad \, \forall (x,y) \in X \times  Y \big \}\]
as well as $\Gamma_{sym}(\mu,\mu)$  the set of symmetric Radon probability measures on
$X \times X$ (i.e., those invariant under the 
permutation  
$R(x_1, x_2)=(x_2, x_1)$) and whose marginals are equal to  $\mu$ on $X.$ We prove the following. 
\begin{theorem}\label{main} Let $X$ and $Y$ be Polish spaces, and let $c : X \times Y \to \R$ be a bounded measurable coupling. Then, for any non-atomic Borel probability 
measure   $ \mu $ on $X$, and any map $T: X \to Y$  such that $(x,z) \to c(x,Tz) $ is upper semi-continuous, we have that 
\begin{equation}
\sup \Big\{ \int_{ X \times X} c(x,Tz) \, d\pi; \, \pi \in \Gamma_{sym}(\mu,\mu)\Big\}=\inf \Big\{ \int_X L\big (x,Tx\big) \, d\mu; L \in {\mathcal L} \Big\}. 
\end{equation}
 Moreover, 
the left-hand side is attained at some transport plan $\pi_0 \in \Gamma_{sym}(\mu,\mu)$, and the right hand-side  
 is attained at some $C$-selfdual Lagrangian $L$, in such a way that 
 \begin{equation}\label{inter}
Tx \in \partial^c_2 H_L(x,z) \qquad \pi_0-a.e. \quad (x,z) \in X \times X.
\end{equation}
If $\pi_0$ is supported on a graph of a measurable map $S:X\to X$, then $S$ is $\mu$-measure preserving, and for $\mu$-almost all $x$ in $X$, we have $S^2x=x$, $H_L(x,Sx)=-H_L(Sx,x)$  and
\begin{equation}\label{inter01}
 Tx \in \partial^c_2 H_L(x,Sx). 
\end{equation}
\item If $T$ is $c$-monotone, then (\ref{inter01}) holds with $S=I$. In this case, $H_L$ is $\mu$-a.e. differentiable in the second variable on the diagonal, and therefore 
$ Tx = \partial^c_2 H_L(x,x)$ or equivalently, 
\begin{equation}
\hbox{$\nabla_2H(x, x)=\nabla_1c(x, Tx)$ \quad for  $\mu$-a.e. $x \in X$.}
\end{equation}
 \end{theorem}
In Section 6, we consider the possibility of using a global variational method to find solutions for equations of the form $p\in Tx$, where $T$ is a given $c$-monotone map. Since $T=\partial_c L(x)$ for some $C$-convex selfdual Lagrangian $L$, the problem reduces to minimizing on $X$ the non-negative functional 
\[
I_p(x)=L\big(x, p\big)-c \big (x, p\big),
\]
 and showing that there exists $x_0$ such that $I_p(x_0)=\inf_{x\in X}I_p(x)=0.$
 For that we needed to make a link with the following notions of arc-wise convexity.
 
Say  that  $F: X \times Y \to \R$ is {\it uniformly  arc-wise convex with respect to the second variable},  if for each $y_0, y_1 \in Y$, there exists a continuous curve $\zeta:[0,1]\to Y$ with $\zeta(0)=y_0$ and $\zeta(1)=y_1$ such that for all $t \in [0,1]$ and all $x \in X$,
\[F(x, \zeta(t))\leq t F(x,y_0)+(1-t)F(x, y_1).
\]
We shall prove the following.

\begin{theorem}\label{Var.0}  Let $X$ be a compact topological space, and let $c: X \times Y \to \R$ be a coupling of $X$ and $Y$, where the latter is a topological space. Suppose $L: X\times Y \to \R$ is a $C$-selfdual Lagrangian and let $H$ be its corresponding antisymmetric Hamiltonian $H$.
\begin{enumerate}
\item  Let $p\in Y$ be such that the function $F_p: X\times X \to \R$ defined by $F_p(x,z):=H(x,z)-c(x,p)$ is lower semi-continuous and uniformly arc-wise convex with respect to the first variable, then the functional  $I_p(x)=L\big(x, p\big)-c \big (x, p\big)$ 
 satisfies  $\inf_{x\in X}I_p(x)=0.$
 
\item If $I_p$ is also lower semi-continuous, then there exists $x_0 \in X$ such that $I_p(x_0)=0$ and $x_0$ is a solution of the equation $p\in \bar \partial_c L(x_0).$
 
\end{enumerate}
\end{theorem}
We note that in the linear case, such a variational principle leads to a resolution of equations  of the form
$Ax\in \partial \phi (x)$,
where $\phi$ is a convex function on a Banach space $X$ and $A: X\to X^*$ is a skew-adjoint operator (i.e., $A^*=-A$). This is done by noticing that $L(x, p)=\phi (x)+\phi^*(Ax +p)$ is a selfdual Lagrangian to which the variational principle B-4 readily applies. However, this reduction  is not possible in a nonlinear setting if we are to solve an equation  of the form  
\begin{equation}\label{equa}
Bx\in \partial \phi_c (x).  
\end{equation}
We are therefore led to solve the equation directly by trying to minimize functionals of the form 
\[
J_p(x)=\phi(x)+\phi^c\big (Bx\big)-c \big (x, Bx\big),
\]
where  $B: X\to Y$ is a {\it $c$-skew adjoint} map, that is if it  satisfies for all $x_1, x_2 \in X$ and  $y_1 \in B x_1 $, $y_2 \in B x_2$, \begin{equation}
\hbox{$c\big (x_1,y_1\big)+c\big ( x_2,y_2\big)=c\big (x_1,y_2\big)+c\big ( x_2,y_1\big)$.}
 \end{equation}
Note that if $c(x,y)=-\|x-y\|^2$ where the norm is given by an inner product, then $B$ is  $c$-skew-symmetric and $B(0)=0$ if and only if $B:X\to X^*$ is a linear skew-symmetric operator, i.e.,
 $\langle Bx,y \rangle=-\langle By,x \rangle$ for all $ x\in X$ and $y \in X^*.$ On the other hand, 
 A simple but non-linear example of a $c$-skew-symmetric map is the counterclockwise rotation by $\pi/2$ on the circle $S^1$, when the cost $c$ is equal the arclength metric $d$. Indeed, note that $d\big(x,Bx\big)=\pi/2$, while 
 $d\big (x,B y \big)+d\big ( y, B x \big)=\pi$ for all $x,y \in S^1.$ \\
 
We shall prove the following.
  
 \begin{theorem}\label{six}  Let $X$ be a compact topological space, $Y$ a topological space, $c: X \times Y \to \R$ a coupling that is uniformly arc-wise convex with respect to the second variable, and $B: X \to Y$  a $c-$skew symmetric map.  If $\phi: X\to \R$ is a lower semi-continuous function such that the map $(x, y)\to \phi(x)-c(x,y)$ is uniformly arc-wise convex in the first variable, then
 \begin{enumerate}
\item $\inf_{x\in X} J_p(x)=0.$ 

\item If $B$ and $c$ are continuous, then the infimum is attained at some  $x_0 \in X$ such that
$Bx_0\in \partial_c \phi(x_0).$
\end{enumerate}
\end{theorem}
Note that since a constant map (i.e., $Bx=p$ for every $x\in X$)  is obviously $c$-skew adjoint, Theorem \ref{six} above could be used to find solutions for equations of the form  $p\in \partial_c \phi(x)$, hence for the inversion of $c$-cyclically monotone operators.

 \section{Metric selfduality up to a transformation}
 Let $U$ and $V$ be two arbitrary sets and $C: U \times V \to \R$ be a finite coupling. The following properties, well known for convex functions, extend easily to this setting. For an arbitrary function $f: U \to \bar \R=\R \cup \{+\infty\}$, the following holds:
 \begin{itemize}
\item $f^{CC}$ is the largest C-convex minorant of $f$.
\item 
(Young inequality) For all $u \in U$ and $v \in V$, we have that $f(u)+f^C(v) \geq C(u,v)$, and 
   $v_0 \in \partial_C f(u_0)$ if and only if  
   $ f(u_0)+ f^C(v_0)= C(u_0,v_0).$
\item If $U$ is an open subset of $ \R^d$, $C$ is differentiable with respect to the first variable, and $f$ is a $C-$convex function on $U$ that is  differentiable at $u_0 \in U$, then 
\begin{equation}
\hbox{$\nabla f(u_0)= \frac{\partial C}{\partial u}(u_0, v_0)$ whenever $v_0 \in \partial_Cf(u_0)$.}
\end{equation}
\item If $g$ is $C-$convex, then
\begin{equation}
\sup_{u \in U} \{g(u)-f(u)\}=\sup_{u \in U} \{g^{CC}(u)-f^{CC}(u)\}.
\end{equation}
\end{itemize}
The following lemma will be used frequently in the sequel. The proof is similar to the proof of Theorem 1.4 in \cite{Moa}.

\begin{lemma}\label{fitz} Let $U$ and $V$ be two arbitrary sets and $C: U \times V \to \R$ be a finite coupling. Assume that $R: V \to U$ is an invertible map from $V$ onto $U$ such that $C(u,v)=C(R^{-1} v, R u).$ If  there exist $\psi: U \to \R \cup \{+ \infty\}$ and $\phi: V \to \R \cup \{+ \infty\}$ such that 
\[C(u,v)\leq \phi(v)+\psi(u) \qquad \forall (u,v) \in U \times V.\]
Then, there exists $L: U \to \R \cup \{+ \infty\}$ such that $L^C=L \circ R,$ and 
\[\hbox{$L(u) \leq \frac{\phi(R^{-1} u)+ \psi(u)}{2}$ for all $u  \in U$ \quad and \quad $L^C(v) \leq \frac{\phi(v)+ \psi(R v)}{2}$ for all $v  \in V.$}
\]
\end{lemma}
\textbf{Proof.} For notational simplicity, write $R_2:=R: V \to U$ and $R_1:R^{-1}: U \to V$ in such a way that $R_2 \circ R_1 =Id_U,$ $R_1 \circ R_2=Id_V$ and $C(u,v)=C(R_1 v, R_2 u).$

 Define $\Phi ( v)=\frac{\phi( v)+\psi(R_2 v)}{2}.$  Denote by $\Phi_C$ the conjugate of $\Phi$ defined by
\[\Phi_C( u)=\sup_{ v \in  V}\{C( u, v)-\Phi( v)\}.\]
It follows from the above together with $R_2 \circ R_1=Id_U$ and $R_1 \circ R_2=Id_V$ that for each $ u \in U$ and $ v \in  V$ we have 
\[\Phi_C( u)+\Phi(R_1 u)+ \Phi_C(R_2 v)+\Phi( v) \geq  C( u, v)+ C(R_2 v, R_1 u)=2 C(u,v).\]
Setting $K( u)= \frac{\Phi_C( u)+\Phi(R_1 u)}{2}$, it follows from the above inequality that
\[K(R_2 v)+K( u) \geq C( u, v), \qquad \text{  for all  }  (u,v) \in U \times V.\]
We also have that $\Phi_C( u)\leq \Phi(R_1 u)$ on $ U.$ In fact,
\begin{eqnarray*}
 C( u, v)-\Phi( v)&=&\frac{C( u, v)-\phi( v)}{2}+\frac{C( u, v)-\psi(R_2 v)}{2}\\
&=&\frac{C( u, v)-\phi( v)}{2}+\frac{C(R_2 v,R_1 u)-\psi(R_2 v)}{2}\\
&\leq&\frac{\psi( u)+\phi(R_1 u)}{2}=\Phi(R_1 u),
\end{eqnarray*}
from which one has $\Phi_C( u)\leq \Phi(R_1 u)$. Then we must have  $ \Phi_C ( u) \leq K( u) \leq \Phi(R_1 u) $ for each $ u \in  U.$
Set
\[\mathcal{W}:= \Big \{W: U \to \R \cup \{+ \infty\};  \quad  \Phi_C ( u) \leq W( u) \leq \Phi(R_1 u) \quad  \& \quad  W(R_2 v)+W( u) \geq C( u, v), \quad \forall (u, v) \in U \times V\Big\}.\]
Note that $\mathcal{W}\not =\emptyset$ as $K \in \mathcal{W}.$ We define an order on the set $\mathcal{W}$ as follows. For $w_1, w_2 \in \mathcal{W} $,
\[w_1 \preceq w_2 \iff w_1(u) \leq w_2(u), \quad \forall u \in U.\]
Note that $(\mathcal{W}, \preceq) $ is a partially  ordered set. We shall use the Zorn's lemma to show that $\mathcal{W}$ has a minimal  element. Let $\{w_i\}_{i\in I}$ be a totally ordered subset of $\mathcal{W}.$ Set $W(u):=\inf_{i\in I} w_i(u)$ and note that  $W\preceq w_i$ for all $i \in I.$ We shall prove that $W  \in \mathcal{W}.$ It is  easily seen that $\Phi_C  \leq W\leq \Phi\circ R_1.$  Let $ u\in U$ and $v \in V.$ For $\epsilon>0,$ there exist $i,j \in I$ such that $W(u) >w_i(u)-\epsilon$ and $W(R_2 v) > w_j(R_2 v)-\epsilon.$ Without loss of generality we may assume that $w_i \preceq w_j.$ It then  follows that
\[W(u)+ W(R_2 v) >w_i(u)-\epsilon +w_j(R_2v)-\epsilon\geq w_i(u) +w_i(R_2v)-2\epsilon\geq C(u,v)-2\epsilon.\]
Since $\epsilon>0$ is arbitrary we obtain  that $W(u)+ W(R_2 v)\geq C(u,v).$ This shows that $W  \in \mathcal{W}.$  Therefore, by the Zorn's lemma  $(\mathcal{W}, \preceq) $ has a minimal element, say  $L.$

Define $\bar L$ to be the conjugate of the function $ v \to L(R_2 v),$ i.e.
\[\bar L( u)=\sup_{ v\in  V} \{C( u, v)-L(R_2 v)\}.\]
We  now  show that $\bar L=L.$   Since $L$ satisfy the inequality $L( u)+L(R_2 v) \geq C( u, v),$ one has $\bar L \leq L.$ On the other hand by virtue of the fact that $\Phi_C(u)\leq L(u)\leq \Phi(R_1 u)$ we obtain 
\begin{eqnarray*}
\bar L( u)=\sup_{ v\in V} \{C( u, v)-L(R_2 v)\} &\geq& \sup_{ v\in V} \{C( u,v)-\Phi( v)\}=\Phi_C( u).
\end{eqnarray*}
Therefore $\Phi_C( u) \leq \bar L( u) \leq L( u) \leq \Phi(R_1 u)$ for all $ u \in U.$ It then follows that $L/2+\bar L/2 \in \mathcal{W}$. Since $L/2+\bar L/2 \leq L$,  the minimality of $L$ yield that $L=\bar L.$  We now show that $L^C=L \circ R_2.$  We have
\begin{eqnarray*}
\bar L(R_2 v)=\sup_{ \tilde v \in V} \{C( R_2 v, \tilde v)-L(R_2 \tilde  v)\}=\sup_{ u\in U} \{C(R_2 v, R_1 u)-L( u)\}=\sup_{ u \in U} \{C( u, v)-L( u)\}=L^C(v),
\end{eqnarray*}
 and therefore $L^C( v) =\bar L(R_2 v)=L(R_2 v)$ on $ V.$
Finally, we have that 
\begin{eqnarray*}
 L( u) 
\leq  \Phi(R_1 u)=\frac{\phi(R_1 u)+\psi( u)}{2},
\end{eqnarray*}
and 
\[L^C(v)=L(R_2 v)\leq \frac{\phi( v)+\psi(R_2 v)}{2}.\]
%\hfill $\square$

 \section{A selfdual representation of $c$-monotone vector fields}
 This section is devoted to the proof of Theorem \ref{three}. To a coupling $c$ on $X\times Y$,  we 
 shall associate the coupling $C$ on the symmetrized space $(X \times Y)\times (Y \times X)$ by the formula
\begin{equation}
C\big((x_1,y_1), (y_2,x_2)\big)=c(x_1,y_2)+c(x_2,y_1).
\end{equation} 
  Here is the main result of this section. 
 
  \begin{theorem}\label{mose} Let $X$ and $Y$ be two sets and $c: X \times Y \to \R$ be  a coupling. \begin{enumerate}
  \item If $L: X \times Y\to \R \cup \{+ \infty\}$ is a $C-$selfdual Lagrangian, then the map $\bar \partial_c L: X \to 2^Y\setminus \{\emptyset\}$ is maximal $c-$monotone.
  
\item Conversely, if  $M$ is a maximal $c-$monotone subset of $X \times Y$,   then there exists a $C-$selfdual function $L: X \times Y \to \R \cup \{+ \infty\}$ such that $M$ is the graph of $\bar \partial_c L$.
\end{enumerate}
\end{theorem}
We first associate to any subset $M$ of $X\times Y$ and any coupling $c: X \times Y \to \R$, a functional that is essentially the counterpart of the Fitzpatrick function in the case of linear coupling. It is defined as follows: $F_{c,M}: X \times Y \to (-\infty, +\infty]$,  defined as  
\begin{equation}
F_{c,M}(x,y):= \sup_{(x_1,y_1) \in M }\Big \{c(x,y_1)+c(x_1,y)-c(x_1,y_1)\Big \}=\sup_{(x_1,y_1)\in M}\Big \{C\big( (x,y), (y_1,x_1) \big)-c(x_1,y_1)\Big\}.
\end{equation}
It is clear that $F_{c,M}$ is $C-$convex. Let us now assume that  $M$ is a maximal  $c-$monotone set.   Then  $c(x,y)\leq F_{c,M}(x,y)$ with equality if and only if $(x,y)\in M.$ It also follows that for all $(x,y)\in X \times Y,$
\begin{eqnarray*}
 F^C_{c,M}(y,x)&=& \sup_{(x_1,y_1) \in X \times Y}\Big \{c(x,y_1)+c(x_1,y)-F_{c,M}(x_1,y_1)\Big \}\\ &\geq & \sup_{(x_1,y_1) \in M}\Big \{c(x,y_1)+c(x_1,y)-c(x_1,y_1)\Big \}=F_{c,M}(x,y).
\end{eqnarray*}
If in addition $(x,y) \in M$ then $c(x,y)= F_{c,M}(x,y)$ and consequently $(y,x)\in \partial_C F_{c,M}(x,y).$ Thus, 
\[ F^C_{c,M}(y,x)+F_{c,M}(x,y)=2c(x,y),\]
from which we obtain $F^C_{c,M}(y,x)=c(x,y).$
Therefore, 
\begin{equation}
c(x,y)\leq F_{c,M}(x,y) \leq F^C_{c,M}(y,x),
\end{equation}
and these inequalities become equalities if and only if $(x,y)\in M.$
The following gives a proof for Part 2) of  Theorem \ref{mose}. 

\begin{proposition}\label{mose.bis} Let $X$ and $Y$ be two sets and $c: X \times Y \to \R$ be a finite coupling. Assume that $M$ is a maximal $c-$monotone subset of $X \times Y$ and that $F_{c,M}$ is the associated Fitzpatrick function to $c$ and $M.$
Then, there exists a $C-$selfdual function $L: X \times Y \to \R \cup \{+ \infty\}$ such that \[c(x,y)\leq F_{c,M}(x,y)\leq  L(x,y) \leq F^C_{c,M}(y,x), \qquad \forall(x,y) \in X \times Y,\] with $ F^C_{c,M}(y,x)=c(x,y)$  if and only if $(x,y)\in M$.
\end{proposition}
\textbf{Proof.} Set $U=X \times Y$ and $V=Y\times X.$ Define $R_1: U \to V$ and $R_2: V\to U$ by $R_1(x,y)=(y,x)$ and $R_2(y,x)=(x,y).$ Note that the symmetrized coupling $C$ associated to $c$ which is defined by 
\[C\big((x_1,y_1), (y_2, x_2) \big)=c(x_1,y_2)+c(x_2,y_1),\]
satisfies the identity $C(u,v)=C(R_2v, R_1u)$ for all $(u,v)\in U \times V.$ Since \[C(u,v)\leq F_{c,M}(u)+F_{c,M}^C(v), \qquad \forall (u,v)\in U \times V,\]
it follows from Lemma \ref{fitz} that there exists a function $L: U \to \R \cup \{+ \infty\}$ with $L^C=L\circ R_2$ such that 
\[L(u)\leq \frac{F_{c,M}(u)+F_{c,M}^C(R_1u)}{2},\qquad  \forall u \in U.\]
This implies that 
\[L(x,y)\leq \frac{F_{c,M}(x,y)+F_{c,M}^C(y,x)}{2} \leq F_{c,M}^C(y,x) \qquad \forall (x,y) \in X \times Y,\]
where we have used that $F_{c,M}(x,y)\leq F_{c,M}^C(y,x)$ for all $x,y.$ It follows from $L(x,y) \leq F_{c,M}^C(y,x)$ that 
\[L^C(y,x)\geq F_{c,M}^{CC}(x,y)=F_{c,M}(x,y).\]
Since $L^C=L \circ R_2$, we obtain that $L(x,y) \geq F_{c,M}(x,y)$ as claimed. \hfill $\square$\\

Now, we establish Part 1 of Theorem \ref{mose}. 
 We isolate the following interesting observation, which connects the $c$-monotonicity of a set in $X \times Y$ to the $C-$cyclical monotonicity of its symmetric enlargement 
 \[E_M:=\Big\{\big ((x,y), (y,x) \big); \, \, (x,y)\in M \Big\},\]
 in the space $(X \times Y) \times (Y \times X)$.
\begin{lemma}\label{enlarg} Let $X$ and $Y$ be two sets and $c: X \times Y \to \R$ be a finite coupling. For a subset $M$ of $X \times Y$, the following assertions hold:
\begin{enumerate}
\item $M$ is $c-$monotone if and only if $E_M$ is $C-$cyclically monotone.
\item $M$ is maximal $c-$monotone in $X\times Y$ if and only if $E_M$ is maximal $C-$cyclically monotone in $E_{X \times Y}$.
\end{enumerate}
\end{lemma}
\textbf{Proof.} Note first that for all $(x_1,y_1)$ and $(x_2,y_2)$ in $X \times Y$ we have
\begin{eqnarray}\label{cuu}
C\big ((x_1,y_1), (y_1,x_1)\big)+ C\big ((x_2,y_2), (y_2,x_2)\big)-C\big ((x_1,y_1), (y_2,x_2)\big)-C\big ((x_2,y_2), (y_1,x_1)\big)=\nonumber\\2\Big( c(x_1,y_1)+c(x_2,y_2)-c(x_1, y_2)-c(x_2, y_1) \Big).
\end{eqnarray}
It follows that if $E_M$ is $C-$cyclically monotone, then it is $C-$monotone and therefore from (\ref{cuu}) we have that $M$ is $c-$monotone.  Now define the function $f: X \times Y \to \R \cup \{+ \infty\}$ by 
\begin{eqnarray*}
	f(x,y)=	\left\{
		\begin{array}{ll}
			c(x,y), & (x,y)\in M \\
			+\infty, & (x,y)\notin M.
				\end{array}
		\right.
	\end{eqnarray*}
If $M$ is $c-$monotone then an easy computation shows that $E_M \subset  Graph \big (\partial_C f\big).$  It then follows that $E_M$ is $C-$cyclically monotone, which proves 1). 

For 2), first assume that $E_M$ is maximal $C-$cyclically monotone. If now $M$ is not maximal $c-$monotone then there exists $(x_1,y_1)\notin M$ such that 
\[c\big (x_1,y_2\big)+c\big ( x_2,y_1\big)\leq c\big (x_1,y_1\big)+c\big ( x_2,y_2\big), \quad \forall  (x_2, y_2) \in M.\]
 It follows from (\ref{cuu}) that  for all $(x_2,y_2) \in M$,
 \[C\big ((x_1,y_1), (y_1,x_1)\big)+ C\big ((x_2,y_2), (y_2,x_2)\big)-C\big ((x_1,y_1), (y_2,x_2)\big)-C\big ((x_2,y_2), (y_1,x_1)\big)\geq 0. \]
 Since $E_M$ is maximal $C-$monotone, then we must have $(x_1,y_1) \in M$ which  leads to a contradiction. 
 
 To prove the other direction we assume that $M$ is maximal $c-$monotone. Let $F_{c,M}$ be the Fitzpatrick function associated to $M$ and $c.$ We shall show that 
 \[E_M=Graph\Big (\partial_C \big(F_{c,M}\big)\Big )\cap E_{X \times Y},\]
 from which the maximal $C-$cyclical monotonicity of $E_M$ in $E_{X \times Y}$ follows.  For that, note that 
 \begin{eqnarray*}
\big( (x,y), (y,x)\big) \in E_M &\iff & (x,y)\in M \\&\iff & F_{c,M}(x,y)=F_{c,M}^C(y,x)=c(x,y)\\
 &\iff & F_{c,M}(x,y)+F_{c,M}^C(y,x)=2c(x,y)\\
 &\iff & (y,x) \in \partial_C \big (F_{c,M}(x,y)\big)\\
 &\iff &\big( (x,y), (y,x)\big) \in Graph\Big (\partial_C \big(F_{c,M}\big)\Big )\cap E_{X \times Y},
  \end{eqnarray*}
as desired. \hfill $\square$\\

\noindent {\bf Proof of Part 1 of Theorem \ref{mose}:} Assume that $L$ is $C$-selfdual, we first see  that $\bar \partial_c L: X \to 2^Y\setminus \{\emptyset\}$ is $c-$monotone. Indeed, let $x_1, x_2 \in X$ and take $y_1 \in \bar \partial_c L(x_1)$ and $y_2 \in \bar \partial_c L(x_2)$. It follows that $L(x_1,y_1)=c(x_1,y_1)$ and $L(x_2,y_2)=c(x_2,y_2).$ This together with the $C-$selfduality of $L$ imply that 
\begin{eqnarray*}
c(x_1,y_1)+c(x_2,y_2)&=&L(x_1,y_1)+L(x_2,y_2)\\
&=&L(x_1,y_1)+L^C(y_2,x_2)\\& \geq & C\big((x_1,y_1), (y_2,x_2)\big)\\
&=& c(x_1,y_2)+c(x_2,y_1),
\end{eqnarray*}
from which the $c-$monotonicity of $\bar \partial_c L$ follows. 

If now $\bar \partial_c L$ is not maximal then there exists a maximal $c-$monotone subset of $X \times Y$ such that $Graph \big(\bar \partial_c L \big) \subset M.$ It follows from Lemma \ref{enlarg} that $E_M$ is maximal $C-$cyclically monotone in $E_{X \times Y}$.
 Let $\partial_C$ be the standard $C-$subdifferential operator and $E_{X \times Y}$  the symmetric enlargement of $X \times Y.$
We shall show that 
\begin{equation}
\label{fz22} Graph \big (\partial_C L \big ) \cap E_{X \times Y} \subseteq E_M.  \end{equation}
Indeed, take $(x,y) \in X \times Y$ such that $(y,x) \in \partial_C L(x,y).$ It implies that 
\[2 c(x,y)=C\big((x,y), (y,x)\big)=L(x,y)+L^C(y,x)=2L(x,y),\]
from which we have that $y \in \bar \partial_c L(x) $ and therefore $(x,y)\in M.$ 
 Therefore, the inclusion (\ref{fz22}) follows. Since $Graph \big (\partial_C L \big ) \cap E_{X \times Y} $ and $E_M$  
are maximal $C-$cyclically monotone in $E_{X \times Y}$ the inclusion in (\ref{fz22}) becomes an equality and therefore $M=Graph\big (\bar \partial_c L\big )$ from which part 1) of Theorem \ref{mose} follows. \hfill $\square$\\

For any Lagrangian $L$ on $X\times Y$, one can associate its Hamiltonian, denoted by $H_L,$ which  is the function on $X \times X$ defined by
 \begin{equation} 
 H_L(z,x)=\sup\{ c(x, y)-L(z,y); \, y\in Y\}.
 \end{equation}
It follows from the definition of $H_L$   that
 \begin{equation}\label{Hamilton1} L^C(y,x)= \sup\{c(z,y)+ H_L(z,x); z \in X\}.\end{equation}
 We also introduce another, possibly set-valued, map $\tilde  \partial H_L $ from $X\to 2^Y\setminus \{\emptyset\}$ by
\begin{eqnarray*}
	\tilde  \partial H_L(x)=	\left\{
		\begin{array}{ll}
			\partial^2_c H_L(x,x), & x \in D_{c,L} \\
			\emptyset, & x \notin D_{c,L}.
				\end{array}
		\right.
	\end{eqnarray*}
The Hamiltonians corresponding to $C$-selfdual Lagrangians have some special features that we list below.

\begin{lemma}\label{Pro-Ham} Let $L: X \times Y \to \R \cup \{+ \infty\}$ be a $C-$selfdual Lagrangian. Then,  the corresponding  Hamiltonian $H_L$ enjoys the following properties:
\begin{enumerate}
\item $x \to H_L(z,x)$ is $c-$convex  for each $z \in X.$
\item $\big[-H_L(.,x)\big ]^{cc}= H_L(x,.)$ for each $x \in X.$
\item $H_L(x,z) \leq -H_L(z,x)$ for all $x,z \in X.$
\item  $L^C(y,x)= \sup\{c(z,y)-H_L(x,z); z \in X\}.$
\item $H_L(x,x)=0$ for all $x \in D_{c,L}.$
\item  The map $\tilde  \partial H_L: X\to 2^Y\setminus \{\emptyset\}$ is maximal $c-$monotone and 
 $\bar \partial_c L=\tilde  \partial H_L.$
 \item  $L^C(y,x)= \sup\{c(z,y)-H_L^{as}(x,z); z \in X\}$ where $H_L^{as}$ is the anti-symmetric Hamiltonian given by
$H_L^{as}(x,z)=\frac{H_L(x,z)-H_L(z,x)}{2}.$

\end{enumerate}
\end{lemma}
\textbf{Proof.}
  (1) follows from the definition. For (2) we fix
$x_1 \in X$ and define $f_{x_1}: X \to \R \cup \{+ \infty\}$ by $f_{x_1}(x)=-H_L(x,x_1).$ It follows that
%\begin{eqnarray*}
\[ f^c_{x_1}(y)= \sup_{x \in X}\{c(x,y)-f_{x_1}(x)\}
= \sup_{x \in X}\{c(x,y)+H_L(x,x_1)\}=L^C(y,x_1)=L(x_1,y),\]
%\end{eqnarray*}
from which one has
\[ f^{cc}_{x_1}(x)= \sup_{y \in Y}\{c(x,y)-f^c_{x_1}(y)\}
= \sup_{y \in Y}\{c(x,y)-L(x_1, y)\}=H_L(x_1,x).
\]
This completes the proof of (2).\\
(3) follows from  (2) together with the fact that $\big[-H_L(.,x)\big ]^{cc}$ is the largest c-convex minorant of $-H_L(.,x).$\\
(4) We have
\begin{eqnarray*}
L^C(y,x)&=& \sup\{c(z,y)+H_L(z,x); z \in X\}\\
&=& \sup\{c(z,y)-\big[ -H_L(z,x) \big]; z \in X\}\\
&=& \sup\{c(z,y)-\big[ -H_L(z,x) \big]^{cc}; z \in X\}\\  
&=& \sup\{c(z,y)-H_L(x,z); z \in X\}.
\end{eqnarray*}
(5) follows from the fact that $H_L(x,z) \leq H_L^{as}(x,z) \leq -H_L(z,x)$ for all $x,z \in X$.\\
To prove 6) consider  $x \in D_{c,L}$, there exist $y_1 \in Y$ such that $L(x,y_1)=c(x,y_1).$  It follows from the definition of $H_L$  
that 
\[H_L(x,x)\geq c(x, y_1)-L(x,y_1)=0.\]
Thus, $H_L(x,x)\geq 0.$ On the other hand by part 3) we have that $H_L(x,x)\leq 0 $
 from which 6) follows.\\
7) It  follows from the 4) above and the $C-$selfduality of $L$ that 
\begin{equation}\label{fz32}L(x,y)=L^C(y,x)=\sup\{c(z,x)-H_L(x,z);\, z \in X\}.
\end{equation}
It follows that for each $x \in D_{c,L}$
\[y \in \partial_c^2 H(x,x) \iff L(x,y)=c(x,y)-H(x,x) \iff L(x,y)=c(x,y) \iff y \in \bar \partial_c L(x),\]
as desired. \hfill $\square$

\begin{corollary} \label{cor1} If  $T$ is a maximal $c$-monotone map from $X$ to $2^Y\setminus \{\emptyset\}$, then there exists a $C$-selfdual function $L$ such that
\begin{equation}\label{inter1}
 Tx = \partial^c_2 H_L(x,x)\qquad  \quad x \in {\rm Dom}(T),
\end{equation}
\end{corollary}

\section{$c-$monotone maps are generically single-valued}

This section is devoted to study   differentiability properties of  sub-antisymmetric functions.  We shall first show that, under certain assumptions on  a smooth manifold $X$,  a polish space $Y$,  a cost function $c: X \times Y \to \R $ and a sub-antisymmetric function $H: X \times X \to \R,$ the map  $\partial^c_2 H$ is single-valued on the diagonal. We then improve a result by Champion-DePascale \cite{de-ch} by showing that if a transport plan is concentrated on a $c$-monotone set, then it is  concentrated on a graph of a  measurable function.    We conclude the section  by proving that  the single-valuedness of  $\partial^c_2 H$ implies the differentiability of $H$ with respect to the second variable on the diagonal.  
Here is the main result of this section.

 \begin{theorem}\label{ap00} Let  $X$ be  a second countable $C^1$ manifold of dimension $d$ with volume measure $\mu$, and let $Y$ be a Polish space. Assume $c: X\times Y \to \R$ is a  measurable coupling that is differentiable with respect to the first variable, while satisfying conditions (\ref{cost1}) and (\ref{cost2}). Assume $H: X\times X  \to \R$ to be a measurable function that is $c-$convex with respect to the second variable such that: \begin{enumerate}
\item  $H$ is sub-antisymmetric, i.e. $H(x,z)+H(z,x)\leq 0$ on $X \times X$ and $H(x,x)=0$ on the diagonal.
\item  There exists a measurable subset $X_0$ of $X$ such that  $\partial^c_2 H(x,x)\not=\emptyset$ for all $x \in X_0.$
\end{enumerate}
Then, $\partial^c_2 H(x,x)$ is single-valued for $\mu$-almost every $x \in X_0$.
  \end{theorem}
 
Here is a direct consequence of the above theorem. 
  \begin{corollary}\label{sing}
 Let $X$ be  a second countable $C^1$ manifold of dimension $d$ with volume measure $\mu$, and let $Y$ be a polish space. Assume $c: X \times Y \to \R$
 is a  measurable coupling that is differentiable with respect to the first variable.
If  $c$ satisfies conditions (\ref{cost1}) and (\ref{cost2}),  then any graph measurable  $c-$monotone map $T: Dom(u) \subset  X \to 2^Y\setminus \{\emptyset\}$ is single-valued on its domain  up to a $\mu$-null set. 
 \end{corollary}
 \textbf{Proof.} It follows from Theorem \ref{three} that there exists a sub-antisymmetric  function $H: X \times X \to \R$ which is zero on the diagonal and $c-$convex with respect to the second variable such that $Tx \in \partial_2^cH(x,x)$ for all $x \in Dom(u).$ It now follows from Theorem  \ref{ap00} that $\partial_2^cH(x,x)$ is single-valued a.e. from which the desired result follows. \hfill $\square$

The following result  
is also an immediate consequence of  Corollary \ref{sing}.
\begin{corollary}\label{champion} {\rm (Champion-De Pascale  \cite{de-ch}):} Let $X$ be a second countable $C^1$ manifold of dimension $d$ and let $Y$ be a
Polish space, $\mu$ and $\nu$ be Borel probabilities respectively over $X$ and $Y$ and let the cost $c$ be
 differentiable  with respect to the first  variable.   We assume that $\mu$ is absolutely continuous with respect to the volume measure  and that $c$ satisfies conditions (\ref{cost1}) and (\ref{cost2}). If a transport plan $\gamma \in \Gamma(\mu, \nu)$ is concentrated on a Borel measurable  $c-$monotone set, then
$\gamma$ is concentrated on a Borel graph.
\end{corollary}

 To prepare for the proof of Theorem \ref{ap00}, we recall the statement of the Yankov-von Neumann-Aumann selection theorem. For a proof, we refer to Hu-Papageorgiou (\cite{HP}, p.158-159). 

\begin{theorem} \label{selec}If $(X, \Sigma)$ is a complete measurable space, $Y$ is a Souslin space,
and $F: X \to 2^Y\setminus \{\emptyset\}$ is graph-measurable, then there exists a sequence $\{f_m\}_{m\geq 1}$ of
$\Sigma-$measurable selectors of $F$ such that
\[ F(x) \subseteq \overline{\{f_m(x)\}}, \qquad m\geq 1,\,\, \forall x \in X.\]
\end{theorem}

\begin{remark} If $\Sigma$ is not complete in the above theorem, then the selectors are  universally measurable.  So if
$(X , \Sigma, \mu)$ is a $\sigma-$finite measure space, we can find a sequence $\{f_m\}_{m\geq 1}$ of $\Sigma-$measurable
selectors of $F$ such that 
\[ F(x) \subseteq \overline{\{f_m(x)\}}, \qquad m\geq 1,\,\, \mu-a.e. \, x \in X.\]
\end{remark}
 \textbf{Proof of Theorem \ref{ap00}.} We may assume without loss of generality,  that $X_0=X,$ as  otherwise, one can perform the same argument on the measure space $(X_0, \Sigma)$ where $\Sigma$ is the restriction of the Borel $\sigma-$algebra of $X$ to $X_0.$ We first assume that  $X$ is   an open set in $\R^d.$ Let $T: X \to Y$ be a universally measurable selection of $\partial^c_2 H$ which exists thanks to Theorem \ref{selec}.  Thus,
 $Tx \in \partial^c_2 H(x,x)$ for all $x \in X$ and therefore  
\begin{equation}\label{ap10}
H(x,z)-H(x,x) \geq c(z,u(x))-c(x,u(x)), \qquad \forall x,z \in X.
\end{equation}
Let $v \in \R^d$ with $\|v\|=1$ and let $\{t_n\}$ be a sequence of positive reals  approaching zero.  For any nonnegative function $g \in C^1_c(X),$ we have
 \begin{eqnarray} \label{exsi1}
 \liminf_{n\to \infty}\int_X g(x)\frac{H(x, x+t_nv)-H(x,x)}{t_n}\,dx &\geq& \liminf_{n\to \infty}\int_X g (x) \frac{c(x+t_nv,u(x))-c(x,u(x))}{t_n}\,dx \nonumber\\
 &\geq& \int_X \liminf_{n\to \infty} g(x) \frac{c(x+t_nv,u(x))-c(x,u(x))}{t_n}dx \nonumber \\&=&\int_X g(x)D_1c(x,u(x))v \, dx,
 \end{eqnarray}
 where we have used Fatou's Lemma,  taking into consideration condition (\ref{cost2})  on the cost $c$. By a simple change of variables, we have for $t$ small enough, 
 \begin{equation}\label{exsi3}
 \int_X g(x)\frac{H(x, x+t_nv)-H(x,x)}{t_n}\,dx=\int_X g(x-t_nv)\frac{H(x-t_nv, x)-H(x,x)}{t_n}\,dx.
 \end{equation}
Since $H(x-t_nv,x)\leq -H(x,x-t_nv)$, it follows that 
\begin{eqnarray}\label{exsi2}
 \limsup_{n\to \infty}\int_X g(x-t_nv)\frac{H(x-t_nv, x)-H(x,x)}{t_n}\,dx&\leq& \limsup_{n\to \infty}\int_X g(x-t_nv)\frac{-H(x, x-t_nv)+H(x,x)}{t_n}\,dx \nonumber \\
 &\leq& \limsup_{n\to \infty}\int_X g (x-t_nv) \frac{c(x,u(x))-c(x-t_nv,u(x))}{t_n}\,dx \nonumber \\
 &\leq& \int_X \limsup_{n\to \infty} g(x-t_nv) \frac{c(x,u(x))-c(x-t_nv,u(x))}{t_n}dx\nonumber \\&=&\int_X g(x)D_1c(x,u(x))v \, dx.
 \end{eqnarray}
Therefore, it follows from (\ref{exsi1}), (\ref{exsi2}) and (\ref{exsi3}) that 
\begin{equation}\label{ap100}\lim_{n\to \infty}\int_X g(x)\frac{H(x, x+t_nv)-H(x,x)}{t_n}\, dx=\int_X g(x)D_1c(x,u(x))v \, dx.\end{equation}
Now choose $\{v_k\}_{k=1}^{\infty}$ to  be a countable dense subset of the sphere $\partial B(0,1).$ By Theorem \ref{selec}, there exists  a sequence $\{f_m\}_{m\geq 1}$ of
$\Sigma-$measurable selectors of $\partial^c_2 H$ such that
\begin{equation}\label{ap101}\partial^c_2 H(x,x) \subseteq \overline{\{f_m(x)\}}, \qquad m\geq 1,\,\, \forall x \in X.\end{equation}
From (\ref{ap100}) it follows that for all integers  $k$ and $m,$  
\[\int_X g(x)D_1c(x,f_1(x))v_k \, dx=\lim_{n\to \infty}\int_X g(x)\frac{H(x, x+t_nv_k)-H(x,x)}{t_n}\, dx=\int_X g(x)D_1c(x,f_m(x))v_k \, dx,\]
from which we have 
\[\int_X g(x)D_1c(x,f_m(x))v_k \, dx=\int_X g(x)D_1c(x,f_1(x))v_k \, dx.\]Since $g$ is arbitrary, we obtain that 
\[D_1c(x,f_m(x))v_k=D_1c(x,f_1(x))v_k, \qquad a.e.\, x \in X.\]
Let $A_{k,m}$ be the full measure  subset of $X$ such that $D_1c(x,f_m(x))v_k=D_1c(x,f_1(x))v_k$ holds for all $x \in A_{k,m}.$ Let $A =\cap_{k,m} A_{k,m}.$ It follows that $A $ is a full measure subset of $X$ and
\[D_1c(x,f_m(x))v_k=D_1c(x,f_1(x))v_k, \qquad \forall x \in A.\]
Since $\{v_k\}_{k=1}^{\infty}$ is a dense subset of $\partial B(0,1),$ it follows  that 
\[D_1c(x,f_m(x))v= D_1c(x,f_1(x))v,\]
for all $x \in A$ and $v \in \partial B(0,1).$ It then follows that $D_1c(x,f_m(x))= D_1c(x,f_1(x))$  for all $x \in A.$ Since $c$ satisfies the twist condition we must have $f_m(x)=f_{1}(x)$ for all $x \in A.$ This together with (\ref{ap101}) imply that $\partial^c_2 H(x,x)=f_1(x)$ for all $x \in A.$  This completes the proof for the case where $X$ is an open subset of $\R^n.$

For  the general case where $X$ is a second countable $C^1$ manifold, we consider a $C^1$ atlas
${(O_i,\Phi_i)}_{i\in \mathbb{N}}$ of $X$ and for each $i$ we define the cost $c_i: \Phi (O_i) \times Y \to \R$ by $c_i = c\circ (\Phi_i^{-1}\times Id).$ We also define $H_i: \Phi_i (O_i)\times \Phi_i (O_i) \to \R$ by $H_i(r,s)=H (\Phi_i(r), \Phi_i(s)).$ Now note that if $x \in \Phi (O_i)$ and  $y \in \partial_2^c H(x,x)$ then $y \in \partial_2^c H_i(r,r)$ where $\Phi_i(r)=x.$
By the previous case $\partial_2^c H_i$ is single-valued a.e. and therefore $\partial_2^c H(x,x)$ is single-valued a.e. with respect to the volume measure on $\Phi_i (O_i)$. 
 \hfill $\square$\\

In the following proposition, we show that  the single-valuedness of $\partial_2^c H$ at a point $x$  gives differentiability of $H$ with respect to the second variable  at $(x,x).$ 

\begin{proposition} \label{diffH}
Let  $X$ be  a second countable $C^1$ manifold, $Y$ be a compact Polish space, and  $c: X\times Y \to \R$ be a   measurable function that is continuously  differentiable with respect to the first variable. Let $H: X\times X  \to \R$ be a locally Lipschitz  function that is $c-$convex with respect to the second variable and such that  its $c-$conjugate with respect to the second variable is continuous.  If $\partial^c_2 H(x,x)$ is single-valued for some  $x \in X$, then $H$ is differentiable with respect to the first variable at the point $(x,x)$, i.e., $D_2 H(x,x)$ exists. \\
 Moreover, by denoting
  \[Dom\big(D_2 H \big):=\{x \in X; \, D_2 H(x,x) \text{ exists}\},\]
  we have that the function  $D_2 H: Dom\big(D_2 H \big) \subseteq X \to T^* X$ is continuous at $x.$ 
 \end{proposition}
 \textbf{Proof.} Since $H$ is locally Lipschitz, Rademacher's theorem yields that $H$ is differentiable almost everywhere on $X \times X$ with respect to the volume measure. Let $Dom\big(D H \big)$ denote the subset of $X \times X$ on which $H$ is differentiable.
 We denote by $D^*_2H(x,x)$ (resp., $\partial_2 H(x,x)$)  the limiting (resp., generalized) Clarke gradients \cite{Cl} with respect to the second variable of $H$ at $x,$ that is 
 \[D^*_2H(x,x)=\big \{\lim_{k\to \infty} p_k;\, p_k=D_2 H(z_k, x_k ), \, (z_k, x_k) \to (x,x),\, (x_k, y_k) \in Dom\big(D H \big)\big \},\]
 and
 \[\partial_2 H(x,x)=conv \big ( D^*_2H(x)  \big).\]
Since $H$ is locally Lipschitz and $Y$ is compact, the sets $D^*_2H(x,x)$ and $\partial_2 H(x,x)$ are non-empty and compact. In order to show that $\partial_2 H(x,x)$ is a singleton, we argue by contradiction and assume that it is not. This implies that $D^*_2H(x,x)$ is not a singleton either and there exist $p,q \in D^*_2H(x,x)$ with $p \not=q.$ Thus there are two sequences $\{(z_k, x_k)\}$ and $\{(z'_k, x'_k)\}$ converging to $(x,x)$ such that $H$ is differentiable at $(z_k, x_k)$ and $(z'_k, x'_k)$ and 
 \[\lim_{k \to \infty} D_2 H(z_k, x_k)=p, \qquad \lim_{k \to \infty} D_2 H(z'_k, x'_k)=q.\]
Let  $L$ be  the $c-$conjugate of $H$ with respect to the second variable, i.e.,
 \[L_H(z_0,y_0)=\sup_{x_0\in X}\{c(x_0, y_0)-H(z_0,x_0)\},\qquad \forall (z_0,y_0)\in X \times Y.\] 
Since $H, L_H$ and $c$ are continuous and $Y$ is compact, we have that $\partial^c_2 H(z_k,x_k)$ and $\partial^c_2 H(z'_k,x'_k)$ are non-empty. Thus, if $y_k \in \partial^c_2 H(z_k,x_k)$ and $y'_k \in \partial^c_2 H(z'_k,x'_k)$
 we must have 
 \[L_H(z_k, y_k)+H(z_k, x_k)=c(x_k,y_k),\]
 and 
 \[L_H(z'_k, y'_k)+H(z'_k, x'_k)=c(x'_k,y'_k).\]
 It then follows that  
 \begin{equation}\label{clark}D_2H(z_k, x_k)= D_1 c(x_k,y_k) \quad {\rm and}\quad  D_2H(z'_k, x'_k)= D_1 c(x'_k,y'_k). \end{equation}
 Again, by the compactness of $Y$ and the continuity of $H, L$  and $c,$ we may assume that the sequences $\{y_k\}$ and $\{y_k'\}$ converge respectively to $y\in \partial_2^c H(x,x)$ and $y'\in \partial_2^c H(x,x).$ Since $\partial_2^c H(x,x)$ is a singleton we obtain that $y=y'.$ By letting $k \to \infty$ in (\ref{clark}), we obtain that  
 \[p= D_1 c(x,y),\qquad q= D_1 c(x,y),\]
and hence $p=q$. This leads to a contradiction and consequently our claim  follows. \hfill $\square$\\

The following result provides  a representation for $c$-monotone maps  an Riemannian manifolds. 
\begin{theorem}\label{cor2} Let $(M,g)$ be a connected compact $C^3$-smooth Riemannian manifold without boundary, equipped with a Riemannian distance $d(x,y)$.  Let $c(x,y)=-d^2(x,y)/2$ and assume that $\mu$ is a Borel probability measure on $M$ that is absolutely continuous with respect to the volume measure.   Then a map $T: X \to X
$ is $c$-monotone $\mu-$a.e.  if and only if there exists an sub-antisymmetric Hamiltonian $H$ that is $c$ convex in the second variable such that  $Tx=\exp_x[\nabla_2 H(x,x)]$ for $\mu$-almost every   $x \in M.$
\end{theorem}
We shall need the following two lemmas. The first is due to McCann \cite{Mn}.   
\begin{lemma}\label{MN}
Let $(M,g)$ be a connected compact $C^3$-smooth Riemannian manifold without boundary, equipped with a Riemannian distance $d(x,y)$. Suppose $\psi=
\psi^{cc}$ for $c(x,y)=-d^2(x,y)/2.$ Then
\begin{equation}\label{MN1}
\hbox{$c(x,y)\leq \psi(x) +\psi^c(y)$ for all $x,y \in M,$}
\end{equation}
 and if $\psi$ is differentiable at a point $x\in M$, then equality holds in $(\ref{MN1})$ if and only if $y=\exp_x[\nabla \psi(x)].$
\end{lemma}
The next lemma  addresses the Lipschitz continuity of  $L$ and $H_L$ required for the application of Proposition \ref{diffH}.
\begin{lemma}\label{lipz} Let $(X,d)$ be a metric space whose diameter $|X|=\sup\{d(x,y); \, x,y \in X\}$ is finite.  Let $c(x,y)=-d^2(x,y)/2$ and assume that $L$ is $C$-selfdual and $H_L$ is its Hamiltonian. Then both $L$ and $H_L$ are Lipschitz continuous.
\end{lemma}
\textbf{Proof.} We first recall the following inequality from (\cite{Mn}, Lemma 1) that
\begin{equation}\label{lipz1}
|d^2(x,y)-d^2(z,y)|\leq 2|X| d(x,z), \qquad \forall x,y,z \in X.
\end{equation}
Since $C \big ((x_1, x_2), (y_1,y_2)\big)=c(x_1,y_1)+c(y_2,x_2)$, 
\[L^C(\tilde v)=\sup \{ C(\tilde x,\tilde v) -L(\tilde x);  \tilde x=(x_1,x_2) \in X \times X\},\]
% it follows from the facts that   $-|X|^2\leq C(\tilde x, v)\leq 0$ is bounded 
and  $L$ is $C$-selfdual,  that $L$ is bounded.  Given 
 $\epsilon >0$ and $\tilde x, v \in X \times X,$ there exists $\tilde z \in  X \times X $  such that $L^C (\tilde v)-\epsilon \leq C(\tilde z,\tilde v) -L(\tilde z) $ and $L^C (\tilde x)\geq C(\tilde z,\tilde x) -L(\tilde z).$
 Thus,
 \[ L^C (\tilde v)-L^C (\tilde x)\leq C(\tilde z,\tilde v)-C(\tilde z,\tilde x)+\epsilon \leq 2|X|d(x_1,y_1)+2|X| d (x_2,y_2)+\epsilon, \]
 where the second inequality follows from (\ref{lipz1}). 
  Since the latter inequality holds for all $\epsilon >0$, the result follows. The Lipschitz property of $H_L$
  follows by a similar argument. \hfill $\square$\\

\noindent \textbf{Proof of Theorem \ref{cor2}.}
If $T$ is $c$-monotone for the coupling  $c(x,z)=-d^2(x,z)/2$, we get from Corollary \ref{cor1} that
\[Tx \in \partial_2^cH_L(x,x), \qquad \mu-a.e. \,\, x \in X,\]
for some $C$-selfdual function $L$ where  $H_L(x,x)=0$ for $\mu$-a.e. $x\in X.$ It then  follows   that
\[L\big (x, Tx \big)+H_L(x,x)=c\big(x , Tx \big) \quad {\rm for}\,\, \mu-a.e. \,\, x \in X.\]
On the other hand for every $x,z \in X$ we have
\[L\big (x, Tx \big)+H_L(x,z)\geq c\big(z , Tx \big).\]
It then follows from Lemma \ref{MN} that for every $x \in X$ where $\nabla_2H_L(x,x)$ exists we must have
$Tx=\exp_x[\nabla_2 H_L(x,x)].$ On the other hand it follows from Lemmas  \ref{lipz}, Theorem \ref{ap00} and Proposition \ref{diffH} that $\nabla_2 H_L(x,x)$ exists $\mu$-a.e. and therefore $Tx=\exp_x[\nabla_2 H_L(x,x)]$ for $\mu$-a.e. $x \in X.$\\
Conversely, if $Tx=\exp_x[\nabla_2 H_L(x,x)]$ for $\mu$-a.e. $x \in X$, we claim that $T$ is then $c$-monotone $\mu$ a.e. Indeed, it follows from Lemma \ref{MN} that 
\[L\big (x, Tx \big)+H_L(x,x)=c\big(x , Tx \big) \quad \mu-a.e. \,\, x \in X.\]
On the other hand, we have  for every $x,y \in X$ that
\[L\big (x, Tx \big)\geq c\big(z , Tx \big)-H_L(x,z).\]
This implies that for $\mu$- a.e. $x, z\in X$, 
\begin{eqnarray*}
c\big(x , Tx \big)+c\big(z , Tz \big)&=&L\big (x, Tx \big)+H_L(x,x)+L\big (z, Tz \big)+H_L(z,z)\\
&=&L\big (x, Tx \big)+L\big (z, Tz \big)\\& \geq &c\big(z , Tx \big)-H_L(x,z)+c\big(x , Tz \big)-H_L(z,x)\\
& \geq &c\big(z , Tx \big)+c\big(x , Tz \big).
\end{eqnarray*}
Thus, 
\[c\big(x, Tx \big)+c\big(y , Tz \big) \geq c\big(z, Tx \big)+c\big(x, Tz \big) \qquad  \mu-a.e.\, \, x,z \in X, \]
 from which the desired result follows.\hfill $\square$

 \section{Generating $c$-monotone fields via symmetric mass transport}

In this section, we show that a symmetric version of the Monge-Kantorovich  transport problem provides a natural way to associate to any map,  in a certain optimal way, a corresponding $c$-monotone re-arrangement. 
If $X$ is a Polish space, we shall denote by  
$\Gamma_{sym}(\mu,\mu)$  the set of Radon probability measures on
$X \times X$ whose marginals are equal to the same probability measure $\mu$ on $X$, and which are invariant under the
cyclic permutation  $R:X \times X \to X \times X$ given by $R(x,z)=(z,x)$. 

If $Y$ is another Polish space, $c$ is a coupling on $X\times Y$, and $C$ is its symmetrized on $(X\times Y)\times (Y\times X)$, we shall denote by ${\mathcal L}$ the class of $C$-selfdual Lagrangians on $X \times Y,$ i.e.
\[{\mathcal L}=\big \{L: X \times Y \to \R \cup \{+ \infty\}; \, L^C(y,x)=L(x,y), \quad \, \forall (x,y) \in X \times  Y \big \}.\]
We now prove the following. 

\begin{theorem}\label{main} Let $(X,d)$ be a metric Polish space, $Y$ another polish space, and  $ \mu $ a non-atomic probability Borel 
measure on $X$. Let $c : X \times Y \to \R$ be a bounded measurable  coupling and $u: X \to Y$  a map  such that $(x,z) \to c(x,u(z)) $ is upper semi-continuous.
Consider the following variational problems:
\begin{eqnarray}
{\rm MK}_{\rm sym}{\rm (c)}:&=&\sup \Big\{ \int_{ X \times X} c(x,u(z)) \, d\pi; \, \pi \in \Gamma_{sym}(\mu,\mu)\Big\}\label{dual}\\
 {\rm DK}_{\rm sym}(\rm c):&=&\inf \Big\{ \int_X L\big (x,u(x)\big) \, d\mu; L \in {\mathcal L} \Big\}.\label{primal}
\end{eqnarray}
The following assertions then hold:
\begin{enumerate}
\item ${\rm MK}_{sym}(c)={\rm DK}_{\rm sym}(\rm c)$ and both of them are attained.

\item If  $\pi_0 \in \Gamma_{sym}(\mu,\mu)$ is a transport plan where  ${\rm MK}_{\rm sym}{\rm (c)}$ is attained, and $L$ is $C$-selfdual Lagrangian where ${\rm DK}_{\rm sym}(\rm c)$ is attained then
 \begin{equation}\label{inter1}
u(x) \in \partial^c_2 H_L(x,z) \qquad \pi_0-a.e. \quad (x,z) \in X \times X.
\end{equation}
 where $\partial^c_2 H_L$ stands for the c-subdifferential of $H_L$ with respect to the second variable. Moreover, if $\pi_0$ is supported on a graph of  a map $S$, then for $\mu$ almost all $x \in X,$ we have
 \begin{equation}\label{inter1}
\hbox{$S^2x=x$,\quad  $H_L(x,Sx)=-H_L(Sx,x)$\quad and \quad 
 $u(x) \in \partial^c_2 H_L(x,Sx)$. }
\end{equation}

\item If $u$ is $c$-monotone, then there exists a $C$-selfdual function $L$ such that $H_L(x,x)=0,$ $\mu$-a.e. $x \in X$, and 
\begin{equation}\label{inter1}
 u(x) \in  \partial^c_2 H_L(x,x)\qquad  \mu-a.e. \quad x \in X.
\end{equation}
\end{enumerate}
\end{theorem}
The proof consists of connecting the above with the standard Monge-Kantorovich theory, which we state in its most general form as established in \cite{V2}. 
\begin{proposition}\label{kan-d}
 Let $(X, \mu)$ and $(Y, \nu)$ be two Polish probability spaces and let  $c: X \times Y \to \R \cup\{-\infty\}$ be an upper semi-continuous cost function such that
\[c(x,y) \leq
a(x)+b(y) \qquad \forall (x,y) \in X \times Y,\]
for some real-valued lower semi-continuous functions
 $ a \in L^1 (\mu), b \in L^1(\nu).$
 \begin{enumerate}
\item The following duality then holds:
\[\max_{\pi \in \Gamma(\mu,\nu)} \int_{X \times Y} c(x,y) \, d \pi(x,y)=\inf_{\phi \in L^1(\mu), \psi  \in L^1(\nu)} \Big (\int_X \phi(x) \, d \mu(x)+\int_X \psi(y) \, d \nu(y)\Big ), \]
where the infimum is over all $\phi$ and $\psi$ such that $c(x,y) \leq \phi(x)+\psi(y)$ on $X \times Y.$
\item  If $c$ and the optimal cost are finite, then there is a measurable $c$-cyclically monotone set $\Gamma \subset X \times Y$ %(which is closed if $a, b$ and $c$ are continuous)
such that any optimal map $\pi \in \Gamma(\mu,\nu)$ is concentrated on $\Gamma.$
\end{enumerate}
\end{proposition}
To prove Theorem \ref{main},  
we shall need a few preliminary results. Let  $X$ and $Y$ be two Polish spaces and $c$ be  a coupling  on $X \times Y.$  Set $U= X \times Y$ and $V=Y \times X.$ We shall again use  the new coupling  $C : U \times V \to \R$  defined by 
\[C\big ((x_1,y_1), (y_2, x_2) \big )=c(x_1,y_2)+c(x_2,y_1).\]
Define  $R_1: U \to V$ and $R_2: V \to U$ by  $R_1(x,y)=(y,x)$ and $R_2(y,x)=(x,y)$ for all $x \in X$ and $y \in Y.$  One can easily deduce that
\[C(u,v)=C(R_2 v,R_1u), \qquad \text{ for all } u=(x_1,y_1)\in U\quad  {\rm and} \quad  v=(y_2,x_2) \in V.\]
 Define  measures $\tilde \mu_1$ and $\tilde \mu_2$ on the Borel $\sigma-$ algebras  $B(U)$ of $U,$ and $B(V)$ of $V$ by
\[\int_{U} f(x,y) \, d \tilde \mu_1= \int_X f\big(x,u(x) \big) \, d \mu,\]
and
\[\int_{V} g(y,x) \, d \tilde \mu_2= \int_X f\big(u(x),x \big) \, d \mu,\]
for all bounded continuous functions $f$ and $g$. It is easily deduced that $R_1\# \tilde \mu_1=\tilde \mu_2$ and $R_2\#\tilde\mu_2=\tilde\mu_1.$  Consider the optimization problem

\[{\rm MK(C, \tilde \mu_1, \tilde \mu_2)}:=\sup \Big \{ \int_{U \times V} C(u,v) \, d\tilde \pi; \, \tilde \pi \in \Gamma(\tilde \mu_1, \tilde \mu_2)\Big\},\]
where $\Gamma(\tilde \mu_1, \tilde \mu_2)$  is the set of Borel probability measures $\tilde \pi$ on $U \times V$ with
$Proj_1(\tilde \pi)= \tilde \mu_1 $ and $ Proj_2( \tilde \pi)= \tilde \mu_2.$\\
It follows from Proposition \ref{kan-d} that ${\rm MK(C, \tilde \mu_1, \tilde \mu_2)}$,
is dual to the following minimization problem 
\[{\rm DK (C, \tilde \mu_1, \tilde \mu_2 )}=\inf \Big\{ \int_{U} \psi(u) \, d\tilde \mu_1 +\int_{V} \phi(v) \, d\tilde \mu_2 ; \, \psi (u)+\phi(v) \geq C(u,v)\Big\}. \]

\begin{lemma}\label{d-til} ${\rm DK (C, \tilde \mu_1,\tilde \mu_2)}$ admits a solution $(L,L^C)$ where $L$ is a $C$-selfdual Lagrangian.
\end{lemma}
\textbf{Proof.} It follows from Proposition \ref{kan-d} that ${\rm DK (C, \tilde \mu_1,\tilde \mu_2)}$ has a solution $(\psi,\phi).$ Thus,
$C(u,v) \leq \phi(v)+\psi(u)$ on $U \times V.$ We also have that  $R_2 \circ R_1=Id_U,$  $R_1 \circ R_2=Id_V$ and $C(u,v)=C(R_2 v, R_1 u)$ for all
 $ u \in U$ and $ v \in  V.$ Thus, it follows from Lemma \ref{fitz} that there exists $L: U\to \R \cup \{+ \infty\}$ with $L^C=L \circ R_2$ such that 
\[C(u,v)\leq L(u)+L^C(v) \leq \frac{\phi(R_1 u)+ \psi(u)}{2}+\frac{\phi(v)+\psi(R_2 v)}{2}, \qquad \forall (u,v) \in U \times V.\]

We now  show that $(L, L^C)$ is a solution of ${\rm DK (C, \tilde \mu_1,\tilde \mu_2)}.$ Indeed, by  integrating the latter inequality on $U \times V$ with respect to the measure $\tilde \mu_1 \otimes \tilde \mu_2$ we get
\begin{eqnarray*}
\int_{U} L( u)i\, d\tilde\mu_1+\int_{V} L^C( u) \, d\tilde\mu_2
\leq \int_{U} \frac{\phi(R_1 u)+\psi( u)}{2} \, d\tilde\mu_1+\int_{ V} \frac{\phi( v)+\psi(R_2 v)}{2} \, d\tilde\mu_2,
\end{eqnarray*}
and since $R_1\#\tilde\mu_1=\tilde\mu_2$ and $R_2\#\tilde\mu_2=\tilde\mu_1$ we have that
\begin{eqnarray*}
\int_{U} L( u) \, d\tilde\mu_1+\int_{V} L^C( v) \, d\tilde\mu_2 \leq \int_{U} \psi( u) \, d\tilde\mu_1+\int_{V} \phi( v) \, d\mu_2.
\end{eqnarray*}
It now follows from the optimality of $(\psi,\phi)$ that the latter is indeed an equality.  \hfill $\square$

\begin{lemma}\label{equ}
Let $A, B$ be two measurable subsets of $X$ and let  $\tilde \pi \in \Gamma(\tilde \mu_1, \tilde \mu_2).$  If $\tilde \pi$ is supported on the graph of a map from $U$ to $Y$ then
\begin{eqnarray}\tilde \pi \big((A\times X)\times (X\times u^{-1}(B)) \big)&=&\tilde \pi \big((A\times X)\times (B \times X) \big) \label{equ1}\\
\tilde \pi \big((u^{-1}(A)\times X)\times (X\times B) \big)&=&\tilde \pi \big((X\times A)\times (X \times B) \big) \label{equ2}
\end{eqnarray}
\end{lemma}
\textbf{Proof.} By assumption $d\tilde \pi(\tilde u,\tilde v)=\delta_{\big (\tilde v=T(\tilde u)\big)}d\tilde \mu_1(\tilde u)$ for some measurable map $T=(T_1,T_2): U \to V$ with the property that $T_\# \tilde \mu_1=\tilde \mu_2.$   Define  measurable maps $F: X \to Y$ and $G: X \to X$ by $F(x)=T_1(x,u(x))$ and $G(x)=T_2(x,u(x)).$ It then follows that
\begin{eqnarray*}
\tilde \pi \big((A\times X)\times (B \times X) \big)&=&\tilde \mu_1 \Big\{\tilde u\in U; \big (\tilde u, T(\tilde u)\big ) \in (A\times X)\times (B \times X)  \Big\}\\
&=& \mu \Big\{ x \in  X; \Big (\big (x, u(x)\big ), T\big (x, u(x)\big )\Big ) \in (A\times X)\times (B \times X)  \Big\}\\
&=& \mu \Big\{ x \in  X; \Big (\big (x, u(x)\big ), \big (F(x), G(x)\big)\Big ) \in (A\times X)\times (B \times X)  \Big\}\\
&=& \mu \Big ( A \cap F^{-1}(B)\Big).
\end{eqnarray*}
By a similar argument we also obtain
\begin{eqnarray}
\tilde \pi \big((A\times X)\times (X\times u^{-1}(B)) \big)=\mu \Big ( A \cap (u\circ G)^{-1}(B)\Big).
\end{eqnarray}
Thus, to prove (\ref{equ1}) we just need to show that $\mu \Big ( A \cap F^{-1}(B)\Big)=\mu \Big ( A \cap (u\circ G)^{-1}(B)\Big).$  To do this we prove that $F=u \circ G.$ Since $T_\# \tilde \mu_1=\tilde \mu_2,$ for every bounded measurable map $f: U \to V$ we have
\[\int_{U} f(T \tilde u) \, d\tilde \mu_1=\int_{V} f(\tilde v) \, d\tilde \mu_2,\]
from which we have
\[\int_{X} f\big(F(x), G(x)\big ) \, d \mu=\int_{X} f\big (u(x),x) \, d\mu.\]
The latter equation for the bounded measurable function \[f(x_1,x_2)=\frac{d\big(x_1,u(x_2)\big )}{1+d\big(x_1,u(x_2)\big )}\]
yields that
\[\int_{X} \frac{d\big(F(x),u(G(x))\big )}{1+d\big(F(x),u(G(x))\big )} \, d \mu=\int_{X} \frac{d\big(u(x),u(x)\big )}{1+d\big(u(x),u(x)\big )} \, d\mu=0.\]
Therefore, $F=u \circ G$ almost surely with respect to the measure $\mu.$ This proves (\ref{equ1}). 

Proof of (\ref{equ2}) is more straightforward than (\ref{equ1}). In fact,
\begin{eqnarray*}
\tilde \pi \big((X\times A)\times (X \times B) \big)
&=& \mu \Big\{ x \in  X; \Big (\big (x, u(x)\big ), \big (F(x), G(x)\big)\Big ) \in (X\times A)\times (X \times B)  \Big\}\\
&=& \mu \Big ( u^{-1}(A) \cap G^{-1}(B)\Big)\\
&=& \mu \Big\{ x \in  X; \Big (\big (x, u(x)\big ), \big(F(x), G(x)\big)\Big ) \in (u^{-1}(A)\times X)\times (X \times B)  \Big\}\\
&=&\tilde \pi \Big(\big (u^{-1}(A)\times X\big )\times \big (X \times B \big ) \Big).
\end{eqnarray*}
This completes the proof.  \hfill $\square$

\begin{lemma}\label{equ3} We have that ${\rm MK}(C, \tilde \mu_1, \tilde \mu_2)=2 {\rm MK}_{\rm sym}(c)$. Moreover, if $\pi_0$ is a maximizer of $ {\rm MK}_{\rm sym}(c)$ then  the plan $\tilde \pi_0$ defined   by \begin{equation} \label{ddd}d \tilde \pi_0\big((x_1,y_1), (y_2,x_2)\big)= \delta_{\big (y_1=u(x_1)\big )} \delta_{\big(y_2=u(x_2)\big )} d \pi_0 (x_1,x_2),\end{equation} is a maximizer of ${\rm MK}(C, \tilde \mu_1, \tilde \mu_2).$
\end{lemma}
\textbf{Proof.} Let $\pi_0$ be a maximizer of $ {\rm MK}_{\rm sym}(c)$ and consider the  plan $\tilde \pi_0$  defined by (\ref{ddd}). It can be easily check that $\tilde \pi_0 \in \Gamma(\tilde \mu_1, \tilde \mu_2).$ It also follows that  
\begin{eqnarray}\label{QQ}{\rm MK}(C, \tilde \mu_1, \tilde \mu_2) &\geq& \int_{U \times V} C\big( (x_1,y_1), (y_2,x_2) \big) \, d \tilde \pi_0 \nonumber\\ &=&\int_{X \times X} c(x_1,u(x_2)) \, d \pi_0(x_1,x_2)
+\int_{X \times  X} c(x_2,u(x_1)) \, d \pi_0(x_1,x_2)=2 {\rm MK}_{\rm sym}(c).\end{eqnarray}
Therefore, $ {\rm MK}(C, \tilde \mu_1, \tilde \mu_2) \geq 2 {\rm MK}_{\rm sym}(c).$\\
 It follows from (\cite{Pr}, Theorem B) that there exists a sequence of transport plans $\{\tilde \pi_n\}_{n \in \N}$, each supported on the  graph  of a measurable map,  such that
 \[\lim_{n \to \infty}\int_{U \times V} C\big( (x_1,y_1), (y_2,x_2) \big) \, d \tilde \pi_n={\rm MK}(C, \tilde \mu_1, \tilde \mu_2).\]
 Define  measures $\pi_n$ on Borel measurable subsets of $X \times X$ by
\[\pi_n (A \times B)= \tilde \pi_n \big ((A \times Y) \times (Y \times B) \big)/2+\tilde \pi_n \big ((B \times Y) \times (Y \times A) \big)/2.\]
 Note that $\pi_n (A \times B)=\pi_n (B \times A)$ and $ Proj_1(\pi_n)= Proj_2(\pi_n)=\mu.$ Therefore $ \pi_n\in \Gamma_{sym}(\mu,\mu).$ By Lemma \ref{equ} we have
 \[\tilde \pi_n \big((A\times X)\times (X\times u^{-1}(B)) \big)=\tilde \pi_n \big((A\times X)\times (B \times X) \big).\]
 and
 \[\tilde \pi_n \big((u^{-1}(A)\times X)\times (X\times B) \big)=\tilde \pi_n \big((X\times A)\times (X \times B) \big).\]
  It then follows that
 \begin{eqnarray*} {\rm MK}_{\rm sym}(c)&\geq& \int_{X \times X } c (x_1,u(x_2)) \, d \pi_n(x_1,x_2)\\
 &=&\frac{1}{2}\int_{U \times V } c (x_1,u(x_2)) \, d \tilde \pi_n\big ((x_1,y_1),(y_2,x_2)\big)+
 \frac{1}{2}\int_{U \times V } c (x_2,u(x_1)) \, d \tilde \pi_n\big ((x_1,y_1),(y_2,x_2)\big)\\
 &=&\frac{1}{2}\int_{U \times V } c (x_1,y_2) \, d \tilde \pi_n\big ((x_1,y_1),(y_2,x_2)\big)+
 \frac{1}{2}\int_{U \times V } c (x_2,y_1) \, d \tilde \pi_n\big ((x_1,y_1),(y_2,x_2)\big)\\
 &=&\frac{1}{2}\int_{U \times V} C\big( (x_1,y_1), (y_2,x_2) \big) \, d \tilde \pi_n \longrightarrow\frac{1}{2} {\rm MK}(C, \tilde \mu_1, \tilde \mu_2)\quad  \qquad (\text{as } n \to \infty).
 \end{eqnarray*}
The latter inequality shows that $2 {\rm MK}_{sym}(c) \geq {\rm MK}(C, \tilde \mu_1, \tilde \mu_2)$ from which together with (\ref{QQ}) we obtain $2 {\rm MK}_{sym}(c) = {\rm MK}(C, \tilde \mu_1, \tilde \mu_2).$ It also  follows from $2 {\rm MK}_{sym}(c) = {\rm MK}(C, \tilde \mu_1, \tilde \mu_2)$ and (\ref{QQ})  that $\tilde \pi_0$ is a maximizer of ${\rm MK}(C, \tilde \mu_1, \tilde \mu_2).$
\hfill $\square$\\

\noindent \textbf{Proof of Theorem \ref{main}.}
Let $\pi_0$ be a maximizer of ${\rm MK}_{sym}(c).$ By  Lemma \ref{equ3}    the plan  $\tilde \pi_0$ defined by  
\[d \tilde \pi_0\big((x_1,y_1), (y_2,x_2)\big)= \delta_{\big (y_1=u(x_1)\big )} \delta_{\big(y_2=u(x_2)\big )} d \pi_0 (x_1,x_2),\] is a maximizer of ${\rm MK}(C, \tilde \mu_1, \tilde \mu_2).$  By Lemma  \ref{d-til}, ${\rm DK (C, \tilde \mu_1, \tilde \mu_2)}$ admits a solution  $(L, L^C)$ where $L$ is a $C$-selfdual Lagrangian. Since ${\rm DK (C, \tilde \mu_1, \tilde \mu_2)}={\rm MK}(C, \tilde \mu_1, \tilde \mu_2)$, one has
 \begin{eqnarray*}
  \int_{U} L(x_1,y_1) \, d \tilde \mu_1+  \int_{V} L^C(y_2,x_2) \, d \tilde \mu_2 &=&\int_{U \times V} C\big( (x_1,y_1), (y_2,x_2) \big) \, d \tilde \pi_0\\
  &=&\int_{U \times V} \big[c(x_1,y_2)+c (x_2,y_1) \big] \, d \tilde \pi_0\\
  &=&\int_{X \times X}\big[ c(x_1,u(x_2))+c (x_2,u(x_1)) \big] \, d\pi_0(x_1,x_2)\\
  &=& 2\int_{X \times X }c(x,u(z)) \, d\pi_0(x,z),
 \end{eqnarray*}
from which we obtain
 \begin{eqnarray*}
  \int_{ X} L \big(x,u(x) \big ) \, d \mu+  \int_{ X} L^C \big (u(x),x \big) \, d \mu =2\int_{X \times X}c(x,u(z)) \, d\pi_0(x,z).
 \end{eqnarray*}
 Since $L^C \big (u(x),x \big)=L \big(x,u(x) \big ) $  we obtain
\begin{eqnarray*}
  \int_{ X} L \big(x,u(x)\big) \, d \mu =\int_{X \times  X} c (x,u(z)) \, d \pi_0.
 \end{eqnarray*}
 We now   show that $\int_{X \times X} H_L(x,z) \, d\pi_0 =0. $
  Since $H_L$ is sub-antisymmetric  and $\pi_0 \in \Gamma^c_{sym}(\mu,\mu)$
\[2\int_{X \times X} H_L(x,z) \, d\pi_0 = \int_{X \times X} H_L(x,z) \, d\pi_0 +\int_{X \times X} H_L(z,x) \, d\pi_0  \leq 0.\]
On the other hand by part (4) of Lemma \ref{Pro-Ham} we have
$L(x,u(x)) \geq c(z,u(x))-H_L(x,z)$ from which we have
\[\int_{X \times X} L(x,u(x)) \, d \pi_0 \geq \int_{X \times X} c(z,u(x))d \pi_0 -\int_{X \times X}H_L(x,z)d \pi_0.\]
Since $\int_{X \times X} L(x,u(x)) \, d \pi_0 = \int_{X \times X} c(z,u(x))d \pi_0$ the above expression implies that $\int_{X \times X}H_L(x,z)d \pi_0 \geq 0$
and therefore the latter is indeed equal to zero i.e. 
\begin{equation}\label{zero} \int_{X \times X}H_L(x,z)d \pi_0 =0.
\end{equation}
 It now follows that
\[\int_{ X} L(x,u(x))  \, d \mu=\int_{X \times X} c(z,u(x))d \pi_0=\int_{X \times X} c(z,u(x))d \pi_0 -\int_{X \times X}H_L(x,z)d \pi_0, \]
and therefore
\[\int_{X \times X} \Big [ c(z,u(x))-H_L(x,z)-L(x,u(x)) \Big ]\, d \pi_0=0. \]
The integrand is non-negative and therefore
\[ c(z,u(x))-H_L(x,z)-L(x,u(x))  =0, \qquad \pi_0- a.e. \quad (x,z) \in X \times X. \]
It then follows  that
\[ u(x) \in \partial_2^c H_L(x,z), \qquad \pi_0-a.e. \quad (x,z) \in X \times X.\]
If now $\pi_0$ is supported on a graph of a map $S$ then
\[0=\int_{X \times X} \Big [ c(z,u(x))-H_L(x,z)-L(x,u(x)) \Big ]\, d \pi_0=\int_{X} \Big [ c(Sx,u(x))-H_L(x,Sx)-L(x,u(x)) \Big ]\, d \mu, \]
and since the integrand  is non-negative one has
 \[u(x) \in \partial_2^c H_L(x,Sx)\qquad  \mu-a.e. \quad x \in X.\]
We now show that $S \in S_2(X).$ Define the  anti-symmetric functional $F$ on $X \times X$ by
\[F(x,z)=\frac{d(x,Sz)}{1+d(x,Sz)}-\frac{d(Sx,z)}{1+d(Sx,z)}.\]
Since $\int F(x,z) \, d\pi_0=0$ one has
\[\int \frac{d(x,S^2x)}{1+d(x,S^2x)} \, d \mu=\int d(Sx,Sx) \, d \mu=0.\]
 This indeed implies that
 $S^2(x)=x$ for $\mu$ almost every $x \in X.$  It follows from (\ref{zero}) that $\int_{X}H_L(x,Sx)d \mu=0.$ Thus,
\[0= 2\int_{X}H_L(x,Sx)d \mu=\int_{X}\big [H_L(x,Sx) +H_L(Sx,x)\big ]d \mu,\] 
and since $H_L(x,Sx) +H_L(Sx,x)\leq 0$ we obtain that $H_L(x,Sx)= -H_L(Sx,x)$ for $\mu$-a.e. $x \in X.$
   This completes the proof.
 \hfill $\square$

3) Suppose now that $u$ is $c$-monotone, then ${\rm MK}_{sym}(c)$  
has a solution  that is supported on the graph of the identity map on $X.$ Indeed, 
since  $u$ is $c$-monotone  we have
\[c\big (x,u(z)\big)+c\big ( z,u(x)\big)\leq c\big (x,u(x)\big)+c\big ( z,u(z)\big), \quad \forall x,z \in X.\]
For every $\pi \in \Gamma_{sym}(\mu,\mu)$ it follows from the latter inequality that
\[\int_{X\times X}c\big (x,u(z)\big)\, d\pi+\int_{X\times X}c\big ( z,u(x)\big)\, d\pi\leq \int_{X\times X}c\big (x,u(x)\big)\, d\pi+\int_{X\times X}c\big ( z,u(z)\big)\, d\pi,\]
from which we obtain
\[\int_{X\times X}c\big (x,u(z)\big)\, d\pi\leq \int_{X\times X}c\big (x,u(x)\big)\, d\mu.\]
This implies that the transport plan $\pi_0$ defined by $d \pi_0(x,z)=\delta_{(z=x)}d\mu(x)$ is a solution of
${\rm MK}_{sym}(c).$  The result then follows from part (2) of Theorem \ref{main}. \hfill $\square$

\section{A variational approach to inverting a $c$-monotone map} 

In this section, we try to extend the variational principle B-4 to a non-linear setting so as to have a global variational method for finding solutions of equations of the form $p\in Tx$, where $T$ is a given $c$-monotone map. Since $T=\partial_c L(x)$ for some $C$-convex selfdual Lagrangian $L$, the problem reduces to minimizing on $X$ the non-negative functional 
\[
I_p(x)=L\big(x, p\big)-c \big (x, p\big),
\]
 and showing that there exists $x_0$ such that $I_p(x_0)=\inf_{x\in X}I_p(x)=0.$
 We shall be able to do so under the following notions of convexity.  

\begin{definition} \rm Let $X$ and $Y$ be topological spaces. 
\begin{itemize}
\item Say that a function $F: X \times Y \to \R$ is {\it arc-wise concave with respect to the second variable}, if for each $y_0, y_1 \in Y$ and any $x \in X$, there exists a continuous curve $\zeta:[0,1]\to Y$ with $\zeta(0)=y_0$ and $\zeta(1)=y_1$ such that for all $t \in [0,1]$, 
\[F(x, \zeta(t))\geq t F(x,y_0)+(1-t)F(x, y_1).\]

\item Say  that  $F: X \times Y \to \R$ is {\it uniformly  arc-wise concave with respect to the second variable},  if for each $y_0, y_1 \in Y$, there exists a continuous curve $\zeta:[0,1]\to Y$ with $\zeta(0)=y_0$ and $\zeta(1)=y_1$ such that for all $t \in [0,1]$ and all $x \in X$,
\[F(x, \zeta(t))\geq t F(x,y_0)+(1-t)F(x, y_1).\]

\item Say also that  $F: X \times Y \to \R$ is {\it geodesically concave} (resp. {\it uniformly geodesically concave}) with respect to the second variable if the curve $\zeta$ can also taken to be a geodesic.
\end{itemize}

Similarly we define arc-wise convex and geodesically convex functions.
\end{definition}
 
\begin{theorem}\label{Var1}  Let $X$ be a compact topological space, $H: X\times X \to \R$ an antisymmetric functional,  and let $c: X \times Y \to \R$ be a coupling of $X$ and $Y$, where the latter is a topological space. Assume that $B: X \to Y$ is a $c-$skew symmetric map such that the functional $F: X\times X \to \R$ defined by $F(x,z):= H(x,z)-c(x,B(z))$ satisfies the following two conditions:
\begin{enumerate}
\item $F$ is uniformly arc-wise convex with respect to the first variable;
\item $F(.,y)$ is lower semi-continuous for each $y \in Y.$
\end{enumerate}
Then,  the functional $I: X \to \R$  defined by 
$I(x)=L_H \big(x, Bx\big)-c \big (x, Bx\big)$ 
 satsifies
\begin{equation} \inf_{x\in X}I(x)=0.
\end{equation}
Moreover, if $I$ is also lower semi-continuous, then there exists $x_0 \in X$ such that $I(x_0)=0$ and $x_0$ is a solution of the following inclusion
\[Bx_0\in \bar \partial_c L_H(x_0).\]
\end{theorem}
We begin by recalling the following topological minimax result, which seems to be suitable to deal with arc-wise convex functions.

\begin{lemma}[K\"onig \cite{Kon}]\label{Ko}  Let $X$ and $Y$ be topological spaces with $X$ compact, and consider $F: X \times Y \to \R\cup \{-\infty, +\infty\}$ such that  for all $y \in Y$, $F(.,y)$ is lower semi-continuous on $X$,  and for all $x \in X$, $F(x,.)$ is upper semi-continuous on $Y$.  Also assume the following conditions:
\begin{enumerate}
\item For any   $\lambda \in \R$, and any nonempty subset $H$ of $X$, the set $\cap_{x\in H}\{y;\, F(x,y)\geq \lambda\}$ is connected. 
\item For any   $\lambda \in \R$, and any nonempty finite subset $K$ of $Y$,  the set $\cap_{y\in K}\{x;\, F(x,y)\leq \lambda\}$ is connected.   \end{enumerate}
 Then, we have
\begin{equation}
\inf_{x \in X}\sup_{y \in Y}F(x,y)=\sup_{y \in Y}\inf_{x \in X}F(x,y). 
\end{equation}
\end{lemma}
Note that hypothesis 1) and 2)  above are satisfied, whenever for all $y \in Y$,  $F(.,y)$ is lower semi-continuous and uniformly arc-wise convex on $X$,  and for all $x \in X$, $F(x,.)$ is upper semi-continuous and uniformly arc-wise concave on $Y$.\\

 \noindent {\bf Proof of Theorem \ref{Var1}:}   Define $G: X \times X \to \R$ by
\[G(x,z):= H(z,x)+c\big(z, B(z)\big)-c\big(x, B(z)\big).\]
By assumptions $(1)$ and $(2)$ the function $G$ is uniformly arc-wise convex with respect to the first variable and $G(.,z)$ is lower semi-continuous for each $z \in Y.$ Since $B$ is $c-$skew symmetric we have 
\[c\big(z, B(z)\big)-c\big(x, B(z)\big)=c\big(z, Bx\big)-c\big(x, Bx\big),\]
and since $H$ is antisymmetric, $G$ can be rewritten as 
\[G(x,z):=-H(z,x)+c\big(z, Bx\big)-c\big(x, Bx\big).\]
This implies that $G$ is uniformly arc-wise concave with respect to the second variable and $G(x,.)$ is upper semi-continuous for each $x \in X.$  It then follows from Lemma \ref{Ko} that 
\[\inf_{x \in X}\sup_{z \in X}G(x,y)=\sup_{z \in X}\inf_{x \in X}G(x,z). \]
On the other hand
\begin{eqnarray*}
\inf_{x \in X}\sup_{z \in X}G(x,z)&=&\inf_{x \in X}\sup_{z \in X}\Big \{-H(x,z)+c\big(z, Bx\big)-c\big(x, Bx\big)\Big\}\\
&=& \inf_{x \in X}\Big\{\sup_{z \in X}\big\{c\big(z, Bx\big)-H(x,z)\big\}-c\big(x, Bx\big)\Big\}\\
&=& \inf_{x\in X} \Big\{L_H\big (x,Bx\big)-c \big (x, Bx\big)\Big\}=\inf_{x\in X} I(x).
\end{eqnarray*}
By a similar argument
\begin{eqnarray*}
\sup_{z \in X}\inf_{x\in X}G(x,z)&=&\sup_{z \in X}\inf_{x \in X}\Big\{-H(x,z)+c\big(y, B(z)\big)-c\big(x, B(z)\big) \Big\}\\
&=&\sup_{z\in X}\Big\{-I(z)\Big\}.
\end{eqnarray*}
Therefore,
\[\inf_{x\in X} I(x)=\sup_{z\in X}\Big\{-I(z)\Big\},\]
from which we obtain $\inf_{x\in X} I(x)=0.$  If now $I$ is lower semi-continuous then there exists $x_0 \in X$ such that $I(x_0)=\inf_{x\in X} I(x)=0.$
It follows that 
\[L_H\big (x_0,Bx_0\big)-c \big (x_0, Bx_0\big)=0,\]
from which we obtain  $Bx_0\in \bar \partial_c L_H(x_0).$ \hfill $\square$

\begin{corollary}\label{Var-c}  Let $X$ and $Y$ be two topological spaces and let $c: X \times Y \to \R$ be a coupling that is uniformly arc-wise convex with respect to the second variable. Assume that $B: X \to Y$ is a $c-$skew symmetric map, and  $\phi: X\to \R$ is a function  such that for each $y\in Y$, the map $x\to \phi(x)-c(x,y)$ is  lower semi-continuous and uniformly arc-wise convex. If $X$ is compact, then 
\[\inf_{x\in X} \big\{\phi(x)+\phi^c\big (Bx\big)-c \big (x, Bx\big)\big\}=0.\]
Moreover, if the infimum is attained at some  $x_0 \in X$, then  $x_0$ is a solution of the following inclusion
\begin{equation}
Bx_0\in \partial_c \phi(x_0).
\end{equation}
\end{corollary}

We also note the following consequence of Lemma \ref{Ko}.  
\begin{proposition} Assume $c: X \times Y \to \R$ is a continuous coupling of two topological spaces $X$ and $Y$ with $X$ compact such that for each $x \in X$, the function $(z,y)\to c(x,y)-c(z,y)$ is uniformly arc-wise concave with respect to the second variable. Let $\phi:X \to \R$ be a lower semi-continuous functional on $X$ such that: 
\begin{enumerate}
\item The function $(x,y)\to \phi(x)-c(x,y)$ is uniformly arc-wise convex with respect to the first variable.
 
\item For every $z,x \in X$, $\sup_{y \in Y} \{c(x,y)-c(z,y)\} \geq \phi(x)-\phi(z).$ 
\end{enumerate}
Then, $\phi$ is $c-$convex. 
\end{proposition}

\textbf{Proof.} For $x \in X$ we have 
\begin{eqnarray*}
\phi^{cc}(x)&=& \sup_{y \in Y}\{c(x,y)-\phi^c(y)\}\\
&=& \sup_{y \in Y}\big\{c(x,y)-\sup_{z \in X}\{c(z,y)-\phi(z)\} \big \}\\
&=& \sup_{y \in Y}\inf_{z \in X}\big\{c(x,y)-c(z,y)+\phi(z)\big \}\\
&=& \sup_{y \in Y}\inf_{z \in X} F_x(z,y),
\end{eqnarray*}
where $F_x(z,y)=c(x,y)-c(z,y)+\phi(z).$ Note that for a fixed $x \in X$ the map $(z,y)\to F_x(z,y)$ satisfies all the assumption of Lemma \ref{Ko} and therefore
\[\inf_{z \in X}\sup_{y \in Y}F_x(z,y)=\sup_{y \in Y}\inf_{z \in X}F_x(z,y). \]
It then follows that 
\begin{eqnarray*}
\phi^{cc}(x)= \sup_{y \in Y}\inf_{z \in X} F_x(z,y)&=&\inf_{z \in X}\sup_{y \in Y}F_x(z,y) \\
&=& \inf_{z \in X}\sup_{y \in Y}\{c(x,y)-c(z,y)+\phi(z)\}.
\end{eqnarray*}
By assumption $(2)$ we have $\sup_{y \in Y}\{c(x,y)-c(z,y)+\phi(z)\}\geq \phi(x)$, from which we obtain
\[\phi^{cc}(x)=\inf_{z \in X}\sup_{y \in Y}\{c(x,y)-c(z,y)+\phi(z)\}\geq \inf_{z\in X}\phi(x)=\phi(x).\]
Therefore, $\phi^{cc}(x)\geq \phi(x)$. Since the inequality $\phi^{cc}\leq \phi$ always holds we indeed have $\phi^{cc}(x)= \phi(x).$  This completes the proof. \hfill $\square$\\

The converse of the latter Proposition  holds under very mild  assumptions.
\begin{proposition} Let $X$ and $Y$ be two sets, $\phi: X \to \R$ any function and $c: X \times Y \to \R$  a coupling. If $c$ is uniformly
arc-wise convex with respect to the second variable, then  the $c$-conjugate $\phi^c$ of $\phi$ is arc-wise convex.  
\end{proposition}

\textbf{Proof.} Take $y_0, y_1 \in Y.$ Since $c$ is uniformly arc-wise convex with respect to the second variable  there exists a continuous curve $\zeta:[0,1]\to Y$ with $\zeta(0)=y_0$ and $\zeta(1)=y_1$ such that 
\[c(x, \zeta(t))\leq t c(x,y_0)+(1-t)c(x, y_1),\qquad \forall t \in [0,1], \, \forall x \in X.\]
It then follows that 
\begin{eqnarray*}
\phi^c(\zeta(t))&=&\sup_{x\in X} \big \{c(x, \zeta(t))-\phi(x) \big \}\\
&\leq & \sup_{x\in X} \big \{tc(x, y_0)+(1-t)c(x,y_1)-\phi(x) \big \}
\\
&\leq & t\sup_{x\in X} \big \{c(x, y_0)-\phi(x) \big \}+(1-t)\sup_{x\in X} \big \{c(x,y_1)-\phi(x) \big \}\\
&=& t\phi^c(y_0)+(1-t)\phi^c(y_1),
\end{eqnarray*}
as desired. \hfill $\square$

\end{document}